\input amstex.tex
\documentstyle{amsppt}
\magnification=\magstep1 \NoBlackBoxes

\pageheight{9.0 true in}

\redefine\phi{\varphi}

\define\ftil{\widetilde{f}}
\define\ta{\theta}
\redefine \phi{\varphi}
\define\R{\Bbb R}
\define\newcommand\koniec{\unskip\nobreak\quad\square\enddemo}

\redefine\epsilon{\varepsilon}
\topmatter
\title Over-rotation numbers for unimodal maps \endtitle


\address   Department of Mathematics,
           University of Alabama in Birmingham,
           University Station,
           Birmingham, AL  35294-2060\endaddress

\email     ablokh\@math.uab.edu  \endemail

\author   A. Blokh and K. Snider \endauthor

\address   Department of Mathematics,
           University of Alabama in Birmingham,
           University Station,
           Birmingham, AL  35294-2060\endaddress

\email  ksnider\@uab.edu \endemail

\thanks One of the authors (A.B.) was partially supported by
NSF grant DMS--0901038. The first version of the paper was produced
when he (A.B.) visited MSRI in the Spring of 1994 (the research at
MSRI was funded in part by NSF grant DMS-9022140).
\endthanks
\abstract{ 
We introduce {\it twist unimodal maps} of the interval and describe
their structure. Sufficient conditions for the growth of
over-rotation interval in families of maps are given.}
\endabstract
\endtopmatter 

\document
\head 0. Introduction\endhead

The present paper is related to a few classical results in
one-dimensional dynamics, namely to a Sharkovski\u\i's theorem on
periods of interval maps and Misiurewicz's results concerning
periodic points of degree one circle maps.

\subhead 0.1 Cycles on the interval and on the circle \endsubhead
To state the\ Sharkovski\u\i\ theorem let us first introduce the
{\it Sharkovski\u\i\ ordering} for positive integers: $$ 3\succ_S
5\succ_S 7\succ_S\dots\succ_S 2\cdot3\succ_S 2\cdot5\succ_S 2
\cdot7\succ_S\dots\succ_S 8\succ_S 4\succ_S 2\succ_S 1\tag{$*$} $$
If $m\succ_S n$, say that $m$ is {\it sharper} than $n$. Let $Sh(k)$
be the set of all integers $m$ with $k\succeq_S m$, let
$Sh(2^\infty)$ be the set $\{1,2,4,\dots\}$, and let $P(\varphi)$ be
the set of (minimal) periods of cycles of a map $\varphi$.

\proclaim{ Theorem S[S]} If $g:[0,1]@>>>[0,1]$ is continuous,
$m\succ_S n$ and $m\in P(g)$ then $n\in P(g)$ and so there exists
$k\in\Bbb N\cup 2^\infty$ with $P(g)=Sh(k)$.  \endproclaim

Similar results for circle maps of degree one are due to Misiurewicz
\cite{M1} who used the notion of the {\it rotation number}. This
notion was first introduced by Poincar\`e \cite{P} for circle
homeomorphisms, then extended to circle maps of degree one by
Newhouse, Palis and Takens \cite{NPT}, and then studied in [BGMY],
[I], [CGT], [M1], [M2], [ALMC], [ALMM] (see [ALM2] with an extensive
list of references).

In fact, one can define rotation numbers in a variety of cases using
the following approach (\cite{MZ}, \cite{Z}). Let $X$ be a compact
metric space, $\phi:X\to\Bbb R$ be a bounded measurable function,
$f:X@>>>X$ be a continuous map. Then for any $x$ the set
$I_{f,\phi}(x)$ of all limits of the sequence ${1\over
n}\sum^{n-1}_{i=0}\phi(f^ix)$ is called the {\it $\phi$-rotation
set} of $x$. Clearly $I_{f,\phi}(x)$ is a closed interval. If
$I_{f,\phi}=\{\rho_\phi(x)\}$ then the number $\rho_\phi(x)$ is
called the {\it $\phi$-rotation number} of $x$; if $x$ is a periodic
point then its rotation number $\rho_\phi(x)$ is well-defined.
Properties of these and related sets in case when $X$ is an interval
are studied in [B2] for a variety of functions $\phi$ and contain much
 information about the dynamics even for an arbitrary function
$\phi$.

For functions $\phi$ related to the dynamics of the map $f$ one
might get additional results about $\phi$-rotation sets; e.g., this
happens for rotation numbers in the circle degree one case [M1]. Let
$f:S^1@>>>S^1$ be a map of degree $1$, $\pi:\Bbb R@>>> S^1$ be the
natural projection which maps the interval $[0,1)$ onto the whole
circle. Fix a lifting $F$ of $f$. Define $\phi_f:S^1@>>> \Bbb R$ so
that $\phi_f(x)=F(X)-X$ for any point $X\in \pi^{-1}x$; then
$\phi_f$ is well-defined, the classical rotation set of a point $z$
is $I_{f,\phi_f}(z)=I_f(z)$ and the classical rotation number of $z$
is $\rho_{f,\phi_f}(z)=\rho(z)$ whenever it exists.

The {\it rotation set} of the map $f$ is $I_f=\cup I_f(x)$; it
follows from [NPT],[I] that $I_f$ is a closed interval (cf. [B2]).
The sum $\sum^{n-1}_{i=0}\phi_f(f^ix)=m$ taken along the orbit of an
$n$-periodic point $x$ is an integer which defines a pair
$(m,n)\equiv rp(x)$ called the {\it rotation pair} of $x$; denote
the set of all rotation pairs of periodic points of $f$ by $RP(f)$.
For real $a\le b$ let $N(a,b)=\{(p,q)\in \Bbb Z^2_+: p/q\in (a,b)\}$
(in particular $N(a,a)=\emptyset$). For $a\in \Bbb R$ and $l\in \Bbb
Z_+\cup \{2^\infty\}$ let $Q(a,l)$ be empty if $a$ is irrational;
otherwise let it be $\{(ks,ns): s\in Sh(l)\}$ where $a=k/n$ with
$k,n$ coprime.

\proclaim{ Theorem M1 [M1]} For a continuous circle map $f$ of
degree $1$ such that $I_f=[a,b]$ there exist $l,r\in \Bbb Z_+\cup
\{2^\infty\}$ such that $RP(f)=N(a,b)\cup Q(a,l)\cup Q(b,r)$.
\endproclaim

The choice of $\phi_f$ is crucial for Theorem M1. By [B1, BM0] an
appropriate choice of $\phi=\phi_f$ leads to results for interval
maps similar to Theorem M1; one can even derive Theorem S from them. More
precisely, let $f:[0,1]\to [0,1]$ be continuous, $Per(f)$ be its set
of periodic points, and $Fix(f)$ be its set of fixed points.  It is
easy to see that if $Per(f)=Fix(f)$ then $\omega(y)$ is a fixed
point for any $y$. Assume from now on that $Per(f)\neq Fix(f)$ and
define a function $\chi_f=\chi$ as follows:

$$\chi(x)=\cases 1/2&\text{if $(f(x)-x)(f^2(x)-f(x))\le 0$,}\\ {0}&\text{if
$(f(x)-x)(f^2(x)-f(x))>0$.}\endcases$$


For any non-fixed periodic point $y$ of period $p(y)$ the integer
$l(y)=$ $\sum^{n-1}_{i=0}\chi(f^iy)$ is at most $p(y)/2$ and is the
same for all points from the orbit of $y$. The pair $orp(y)=(l(y),
p(y))$ is called the {\it over-rotation pair} of $y$, and {\it
coprime over-rotation pair} if $p,q$ are coprime. The set of all
over-rotation pairs of periodic non-fixed points of $f$ is denoted
by $ORP(f)$ and the $\chi$-rotation number
$\rho_\chi(y)=\rho(y)=l(y)/p(y)$ is called the {\it over-rotation
number} of $y$. Observe that by Theorem S and by the assumption that
$Per(f)\neq Fix(f)$ it follows that $f$ has a point of period $2$
and that the over-rotation number of this point is $1/2$; in other
words, the set of all over-rotation numbers of periodic points of
$f$ includes $1/2$ and, therefore, $1/2$ belongs to the union of all
$\chi$-rotation sets $I_{f, \chi}(x)$ defined earlier.

We introduce the partial ordering $\Vdash$ among all pairs of
integers $(s,t),\,0<s\le t/2$: $(p,q)\Vdash (k,l)$ if $k/l\in (p/q,
1/2]$.

\proclaim{Proposition BM1 [BM0]} If $(p,q)\Vdash (k,l)$ and
$(p,q)\in ORP(f)$ then $(k,l)\in ORP(f)$.
\endproclaim

This proposition implies Theorem~S. Indeed, let $f$ be an interval
map and consider odd periods. For any $2n+1$ the closest to $1/2$
over-rotation number of a periodic point of period $2n+1$ is
$\frac{n}{2n+1}$.  Clearly
$\frac{n}{2n+1}<\frac{n+1}{2n+3}<\frac12$.
Hence for any periodic point $x$ of period $2n+1$ its over-rotation
pair
$orp(x)$ is $\Vdash$-stronger than the pair $(n+1,2n+3)$, 
and by the proposition the map $f$ has a point of period $2n+3$.
Also, for any $m$ we have $(n,2n+1)\Vdash (m,2m)$, so by the same
proposition the map $f$ has a point of period $2m$. Applying this
result to the maps $f,f^2,f^4,\dots$ one can prove the
Sharkovski\u\i's theorem for all periods but the powers of $2$;
additional arguments covering the case of powers of $2$ are quite
easy. Below we extend the definition of the $\Vdash$-ordering;
namely, if $0<k\le n/2$ then: \roster
\item  if $(p,q)\Vdash (k,l)$ then $(p,q)\gtrdot (k,l)$,
\item if $p/q=k/l=m/n,\,m,n$ coprime, then $(p,q)\gtrdot
(k,l)$ if and only if $(p/m)$ is sharper than $(k/m)$ (both $(p/m)$
and $(k/m)$ are integers).
\endroster

\proclaim{Theorem BM2 [BM0]}  If $(p,q)\gtrdot (k,l)$ and $(p,q)\in
ORP(f)$ then $(k,l)\in ORP(f)$.
\endproclaim

Theorem~BM2 implies a full description of sets $ORP(f)$ for interval
maps, close to that from Theorem~M1.  It implies that the closure of
the set of over-rotation numbers of periodic points of $f$ is an
interval $I_f=[\rho_f, 1/2], 0\le \rho_f\le 1/2,$ called the {\it
over-rotation interval} of $f$. Observe that over-rotation numbers
of periodic points of $f$ are simply $\chi$-rotation numbers of
those points. It is natural to consider the connection between
$I_f$, defined by the over-rotation numbers of {\it periodic points}
of $f$, and the union of {\it all} $\chi$-rotation sets $I_{f,
\chi}(x)$ defined earlier. Before we describe
this connection in detail, we would like to discuss certain
difficulties related to such description; hopefully, this will help the
reader understand assumptions and exclusions which are necessary
here.

For over-rotation numbers, the dynamics in small neighborhoods of
fixed points can play a misleading role. First of all, a point $x$
which maps into a fixed point yields a sequence $\chi(f^i(x))$ which
eventually consists of zeros and hence yields the set
$I_{f,\chi}(x)=\{0\}$. However, this clearly has no implications for
the dynamics of periodic points of the map. Hence {\it pre-fixed}
points $x$ should not be considered as we compare $I_f$ and the
union of {\it all} $\chi$-rotation sets $I_{f, \chi}(x)$.

In general, points which contain some fixed points in their limit
sets should not be considered here because then the seemingly rich
dynamics may take place over a shrinking to zero part of the space
and therefore should be ignored rather than taken into account. To
explain this, let us draw analogy with the case of the entropy. It
is known that for continuous interval maps it can happen so that the
entropy of such maps is large (even infinite) while it is assumed on
smaller and smaller invariant sets converging to fixed points of the
map. Similarly, it can happen that the dynamics in a small
neighborhood of, say, an attracting fixed point $a$ is chaotic. That
may lead to a rich set of sequences $\chi(f^i(x))$ and large
$\chi$-rotational sets of such points while having absolutely no
bearing on the set of periodic points of the map at all (e.g., maps
like that can be such that {\it all} points are attracted to a
unique fixed point). To avoid this ``artificial'' richness we
consider only {\it admissible} points.

More precisely, by a {\it limit measure} of a point $x$ we mean a
limit of ergodic averages of the $\delta$-measure concentrated at
$x$; clearly, any limit measure is invariant. The results of [B2]
apply to a wide variety of functions, but we only state them as they
apply to the function $\chi_f=\chi$ and over-rotation numbers. Call
a point $x$ {\it admissible} if any limit measure $\mu$ of $x$ is
such that $\mu(Fix(f))=0$; since $\mu$ is invariant, this implies
that in fact the set of all points $x$ which are eventual preimages
of fixed points of $f$ is of zero $\mu$-measure. Since the set of
discontinuities of $\chi$ is contained in the union of the set of
fixed points $Fix(f)$ of $f$ and their preimages, we see that for an
admissible point $x$ the set of discontinuities of $\chi$ is of zero
limit measure for any limit measure of $x$.


The connection between $I_f$ and the union of $\chi$-rotation sets
$I_{f, \chi}(x)$ for all admissible points $x$ is established in
papers [B2, B3] and illustrated in Theorem B1; the part of
Theorem~B1 concerning rational rotation numbers and periodic points
follows from the definitions and Theorem BM2. To state the last part
of Theorem B1 we define {\it piecewise-monotone} interval maps. Say
that $f:I=[a, b]\to \Bbb R$ is a {\it piecewise-monotone} map if
there are (perhaps degenerate) closed intervals $a\le
C_0=C<C_1<\dots<C_l\le b$ with $f$ on each component of $[a,
b]\setminus \bigcup C_i$ being monotone and $C_i, 0\le i\le l$ being
a (perhaps, degenerate) flat spot for each $i$. A flat spot is an
interval I such that f|I is a constant. Sets $C_i$ are said
to be {\it critical} (sets of $f$). A degenerate set $C_i$ is called
a {\it critical point} of $f$.  This determines $l=l(f)$; components
of $[a, b]\setminus \bigcup C_i$, denoted by $I_1, \dots, I_l$, are
called {\it laps} of $f$. Thus, laps and sets $C_0\cup I_1, C_2,
I_2, \dots, C_{l-1}, I_l\cup C_l$ form a partition of $[a, b]$
called a {\it basic} partition of $[a, b]$. For simplicity, we
always assume that $f$ is a piecewise-monotone map with only
finitely many fixed points of $f$ as well as finitely many their
first preimages.

\proclaim {Theorem B1 [B2, B3]} The following statements are true.
\roster
\item If $f$ is continuous and $\rho_f<1/2$ then for any $a\in (\rho_f, 1/2]$
there is an admissible point $x$, generic for a measure $\mu$, such
that $I_f(x)=\{a\}$. Moreover, if $p, q$ are positive integers such
that $p/q\in (\rho_f, 1/2]$ then we may choose $x$ to be periodic
with over-rotation pair $orp(x)=(p,q)$.
\item If $x$ is an admissible point, then $I_f(x)\subset I_f=[\rho_f, 1/2]$.
\item If $f$ is piecewise-monotone and $\rho_f\ne 0$ then there exists
an invariant measure $\mu$ such that $f$ is minimal on the support
of $\mu$ and there exists a point $x$, generic for $\mu$ and such
that $I_{f, \chi}(x)=\{\rho_f\}$.
\endroster
\endproclaim

Theorems BM1 and B1 guarantee the existence of a periodic orbit of
any over-rotation pair $(p,q)$ with rational $p/q\in int\,I_f$. By
Theorem B1 there is also a point whose over-rotation set coincides
with a given irrational number from $int\,I_f$. The remaining case
of the left endpoint of $I_f$ is resolved for piecewise-monotone
maps in Theorem B1(3). It is easy to see that the piecewise-monotone
assumption is necessary here. Indeed, consider the following
example. Let $f:[0, 1]\to [0, 1]$ be a map with a sequence of
concatenated invariant intervals $I_j$ such that $\bigcup I_j=[0,
1)$. Then $1$ is clearly an $f$-fixed point. Suppose that $f|_{I_j}$
is piecewise-monotone for any $j$ and such that the over-rotation
intervals of $f|_{I_j}$ grow to their union $(\rho_f, 1/2]$. Then
there exists no admissible point $x$ for which $I_{f, \chi}(x)$
contains $\rho_f$. Thus in Theorem B1(3) above the
piecewise-monotone assumption is necessary.

In a recent paper by Jozef Bobok [Bo] the case covered in Theorem
B1(3) is studied in great detail and depth resulting into a much
more precise claim. The results of [Bo] which complement and further
develop Theorem B1 are summarized in Theorem Bo1 below. Recall that
a dynamical system is said to be {\it strictly ergodic} if it has a
unique invariant measure. To state Theorem Bo1 in full generality we
need a couple of notions on which we will elaborate in Subsection
0.2.

A {\it pattern} is a cyclic permutation of the set
$T_n=\{1,2,\dots,n\}$. A pattern $\pi$ {\it forces} a pattern
$\theta$ if a continuous interval map $f$ which exhibits $\pi$ also
exhibits $\theta$. By \cite{Ba} forcing is a partial ordering. One
can talk about the {\it over-rotation pair} $orp(\pi)$ and the {\it
over-rotation number} $\rho(\pi)$ of a pattern $\pi$. We call a
pattern $\pi$ an {\it over-twist pattern} (or just an {\it
over-twist}) if it does not force other patterns of the same
over-rotation number.

\proclaim {Theorem Bo1 [Bo]} Let a point $x$ and a measure $\mu$ be
as defined in Theorem {\rm B1(3)}. Then the map $f|_{\omega(x)}$ is
strictly ergodic with $\mu$ being the unique invariant measure of
$f|_{\omega(x)}$. Moreover, if $\rho_f$ is rational then $x$ is
periodic and the pattern of the orbit of $x$ is an over-twist
pattern of over-rotation number $\rho_f$.
\endproclaim


Thus, numerical information about a map, compressed to $I_f$,
implies various types of the limit behavior of points reflected by
their rotation numbers. Can one say more about the dynamics of a map
$f$ if $I_f$ contains some number $a$? This question is addressed in
the present paper for unimodal maps.

\subhead 0.2 Patterns and unimodal maps \endsubhead We need
definitions from {\it combinatorial dynamics} ([ALM2]). A map $f$
has a {\it horseshoe} if there are points $a,b,c$ such that either
$f(c)\le a=f(a)<b<c\le f(b)$ or $f(c)\ge a=f(a)>b>c\ge f(b)$. It is
easy to see ([BM0]) that if a map has a horseshoe then it has
periodic points of all possible over-rotation numbers.
A map (not even necessarily one-to-one) of the set $T_n$ into itself
is called a {\it non-cyclic pattern}. If an interval map $f$ on its
cycle $P$ is conjugate to a pattern $\pi$ by an increasing map then
$P$ is a {\it representative} of $\pi$ in $f$ and $f$ {\it exhibits}
$\pi$ on $P$; if $f$ is {\it monotone (linear)} on each
complementary to $P$ interval, we say that $f$ is {\it $P$-monotone
($P$-linear)}(\cite{MN}).

A pattern $\pi$ is said to have a {\it block structure} if there is
a collection of pairwise disjoint segments $I_0, \dots, I_k$ with
$\pi(T_n\cap I_j)=T_n\cap I_{j+1}, \pi(T_n\cap I_k)=T_n\cap I_0$;
the intersections of $T_n$ with intervals $I_j$ are called {\it
blocks} of $\pi$. A pattern without a block structure is said to be
{\it irreducible}. If we identify blocks, we get a new pattern
$\pi'$, and then $\pi$ is said to have a block structure {\it over
$\pi'$}. A pattern $\pi$ {\it forces} a pattern $\theta$ if a
continuous interval map $f$ which exhibits $\pi$ also exhibits
$\theta$. By \cite{Ba} forcing is a partial ordering. If $\pi$ has a
block structure over a pattern $\ta$, then $\pi$ forces $\ta$. By
[MN] for each pattern $\pi$ there exists an irreducible pattern
$\pi'$ over which $\pi$ has block structure (in particular, $\pi'$
is forced by $\pi$).

One can talk about the {\it over-rotation pair} $orp(\pi)$ and the
{\it over-rotation number} $\rho(\pi)$ of a pattern $\pi$. We call a
pattern $\pi$ an {\it over-twist pattern} (or just an {\it
over-twist}) if it does not force other patterns of the same
over-rotation number. Theorem BM2 and the properties of forcing
imply the existence of over-twist patterns of a given rational
over-rotation number between $0$ and $1$: it implies that a map
which has a periodic point of rational over-rotation number $\rho$
exhibits an over-twist pattern of rotation number $\rho$. By Theorem
BM2 an over-twist pattern has a coprime rotation pair; in
particular, over-twists of rotation number $1/2$ are of period $2$,
so from now on we consider over-twists of over-rotation numbers
distinct from $1/2$. Combining this with Theorem Bo1 and Theorem B1,
we come up with the following way of describing periodic dynamics of
a piecewise monotone interval map $f$: {\it if $f$ has the
over-rotation interval $I_f$ then for any rational number $\rho\in
I_f$ there exists an $f$-periodic point $x$ whose orbit exhibits an
over-twist pattern of over-rotation number $\rho$.} This explains
why studying over-twist patterns is important.

In the setting of (non-cyclic) patterns it is useful to consider an
interpretation of over-rotation numbers which is close to symbolic
dynamics. The following construction is a key ingredient of
one-dimensional combinatorial dynamics. Let $\pi$ be a (non-cyclic)
pattern, $P$ be a finite set with a map $f:P\to P$ of (non-cyclic)
pattern $\pi$ and $f$ be a $P$-linear map. Say that a component $I$
of $[0, 1]\setminus P$ {\it $\pi$-covers} another such component $J$
if $J\subset f(I)$. Construct the oriented graph $G_\pi$ whose
vertices are components of $[0, 1]\setminus P$ and whose edges
(arrows) go from $I$ to $J$ if and only if $I$ $\pi$-covers $J$.
Clearly, $G_\pi$ does not depend on the actual choice of $P$ and the
definition is correct.

A cycle (and the pattern it represents) is {\it divergent\/} if it
has points $x<y$ such that $f(x)<x$ and $f(y)>y$. A cycle (pattern)
that is not divergent will be called {\it convergent\/}. It is
well-known that a pattern does not force a horseshoe if and only if
it is convergent (the main ideas of the proof date back to the
original paper by Sharkovski\u\i~[S]). Suppose that $\pi$ is a
convergent pattern and that $P$ is a periodic orbit of pattern
$\pi$. Let $f$ be a $P$-linear map. Then $f$ has a unique fixed
point $a$. Consider the set $Q=P\cup \{a\}$ and denote its pattern
by $\pi'$. We will work with the oriented graph $G_{\pi'}$.

Suppose that there is a real-valued function $\psi$ defined on
arrows of $G_{\pi'}$. This is a classical situation of
one-dimensional symbolic dynamics. It is well-known [ALM2] then that
the maximal and the minimal averages of $\psi$ along all possible
paths (with growing lengths) in $G_{\pi'}$ are assumed on periodic
sequences. In particular, if the values of $\psi$ on arrows are all
rational, then the maximum and the minimum of those averages are
rational. We choose a specific function $\psi$ as follows. Associate
to each arrow in $G_{\pi'}$ the number $1$ if it corresponds to the
movement of points from the right of $a$ to the left of $a$.
Otherwise associate $0$ to the arrow. As explained above, this
yields rational maximum and rational minimum of limits of averages
of $\psi$ taken along all possible paths (with growing lengths) in
$G_{\pi'}$, and these extrema are assumed on periodic sequences.

Now we define unimodal maps. Let $I=[0,1]$.  A continuous map
$f:I\to I$ is {\it unimodal} if there is $c\in(0,1)$ such that
$f(c)=1, f(1)=0,\, f$ increases on $[0,c]$ and decreases on $[c,1]$
(a unique critical point of a unimodal map $g$ is denoted by
$c_g\equiv c$). Suppose there is a fixed point $b\in [0,c]$. Then
the map has a horseshoe (indeed, $0=f^2c\le fb=b<c<fc=1$) which by
the above implies that $f$ has periodic points of all over-rotation
pairs, $I_f=[0,1/2]$, and all the points to the left of $b$ are
attracted by fixed points. To avoid considering this trivial case,
from now on we assume that there are no fixed points in $[0,c)$.
Given a (unimodal) interval map $f$ we consider its over-rotation
interval $I_f=[\rho_f, 1/2]$. In this paper we study $\rho_f$ for
unimodal maps. By the previous paragraph if $c$ is periodic or
preperiodic, then $\rho_f$ is rational. Moreover, similarly to the
previous paragraph it follows that if $c$ is attracted to a periodic
point of $f$ then still $\rho_f$ is rational. The main idea of the
paper is to study $\rho_f$ by first constructing a discontinuous
lifting $F$ whose over-rotation numbers coincide with those of $f$
and then study $F$ in the spirit of [M3].

Fix a unimodal map $f$. Next we briefly describe the main notions of
{\it kneading theory}, due to Milnor and Thurston [MT]. For each
point $x\in [0, 1]$ we define its {\it itinerary} as the sequence
$i(x)=i_0(x)i_1(x)\dots$ of symbols $L, C$ or $R$ so that $i_j(x)=L$
if $f^j(x)<c$, $i_j(x)=C$ if $f^j(x)=c$, and $i_j(x)=R$ if
$f^j(x)>c$. Define the order of the symbols to be $L < C < R$. Now,
suppose that $A=a_0\dots$ and $B=b_1\dots$ are two sequences of
symbols $L, C$ or $R$. Define the order among them as follows.
Choose the smallest $j$ with $a_j\ne b_j$. Then we set $A\succ B$ if
there is an even number of $R$'s among $a_0, \dots, a_{n-1}$ and
$a_n>b_n$ or if there is an odd number of $R$'s among $a_0, \dots,
a_{n-1}$ and $a_n<b_n$. It is shown in [MT] that $x>y$ implies
$i(x)\succeq i(y)$ and $i(x)\succ i(y)$ implies $x>y$.

An itinerary $A$ is said to be {\it shift maximal} if $A\succeq
\sigma^j(A)$ for any non-negative $j$ where $\sigma$ is the left
shift. The {\it kneading sequence} of $f$ is the itinerary
$\nu(f)=i(f(c))$. Clearly, $\nu(f)$ is shift maximal. By [MT],
$\nu(f)\succeq A\succeq \sigma(\nu(f))$ if and only if there exists
a point $x$ such that $i(x)=A$. Therefore, the over-rotation
interval of a unimodal map $f$ is determined by its kneading
sequence and we can talk about the over-rotation interval $I_\nu$ of
a kneading sequence $\nu$. In fact, by definition and properties of
kneading sequences if $\nu_2\succ \nu_1$ and $f_1, f_2$ are
respective unimodal maps, then all patterns exhibited by
$f_1$-cycles are also exhibited by $f_2$-cycles. In particular, then
$I_{\nu_2}\subset I_{\nu_2}$.

\subhead 0.3 Main results \endsubhead Let us proceed in a more
detailed manner. Consider the shift $\zeta_\rho$ by $\rho$ on
$[0,1)$ modulo $1$ assuming $\rho\le 1/2$. Define the kneading
sequence $\nu_\rho=(\nu_\rho(0), \nu_\rho(1), \dots)$ as follows:
$$\nu_\rho(n)=\cases {C}&\text{if $\zeta_\rho^{n+1}(\rho)=0$,}\\ {R}&\text{if
$0<\zeta_\rho^{n+1}(\rho)<2\rho$}\\ {L}&\text{if $2\rho \le
\zeta_\rho^{n+1}(\rho)$}\endcases$$ Clearly, we have $\nu_\rho=(R,
L, \dots)$ (except for $\nu_{1/2}=(R, C, R, C, \dots)$). We show
that $\nu_\rho$ is a kneading sequence of a unimodal map (this can
be done formally, but we prefer a more geometrical approach). Note
that $\nu_\rho$ is periodic if $\rho$ is rational and non-periodic
otherwise. Let $\rho=p/q,\, p, q$ coprime.  It is easy to see that
the corresponding to $\nu_\rho$ pattern $\gamma_\rho$ is given by
the cyclic permutation $\phi$ of the set of points $\{0,1/q, \dots,
(q-1)/q\}$ where $\phi$ is defined as follows: \roster
\item $\phi(j/q)=j/q+p/q$ for $0\le j\le q-2p-1$,
\item $\phi(j/q)=(2q-2p-1-j)/q$ for $q-2p\le j\le q-p-1$,
\item $\phi(j/q)=(q-1-j)/q$ for $q-p\le j\le q-1$.
\endroster

In other words, here is what the pattern $\gamma_{p/q}$ does with
$q$ points $x_0, \dots, x_{q-1}$ of the periodic orbit. The first
$q-2p$ points from the left are shifted to the right by $p$ points.
The next $p$ points are ``flipped'' (i.e. the orientation is
reversed, but the points which are adjacent remain adjacent) all the
way to the right. Finally, the last $p$ points of the orbit are
flipped all the way to the left. As we will see in what follows {\it
the pattern $\gamma_\rho$ is the only unimodal twist pattern of
rotation number $\rho$} in the case of a rational $\rho$. Note that
the same fact follows from the results of [BK] where different
methods are used.

Let us now introduce a non-cyclic pattern $\gamma'_\mu$ for
$\mu<1/2$. \roster
\item $D_\mu=\{0, 1/q, \dots, \dfrac {q-1}q, a', a\}$ where
$\dfrac {q-2p-1}q<a'<\dfrac {q-2p}q, \dfrac {q-p}q<a<\dfrac
{q-p+1}q$; also let $j=\gamma^{-1}(q-2p)$.
\item $\gamma'_\mu: D_\mu\to D_\mu$.
\item $\gamma'_\mu(i/q)=\gamma_\mu(i/q)$ for any $i/q\neq j/q,
a, a'$.
\item $\gamma'_\mu(j/q)=a', \gamma'_\mu(a')=\gamma'_\mu(a)=a$.
\endroster

As follows from the definition, the pattern $\gamma'_\rho$ can be
easily obtained from the pattern $\gamma_\rho$. Indeed, consider the
pattern $\gamma_\rho$ with added points $a', a$ at the appropriate
places. Change the map on the first preimage of $c$ so that it maps
to $a'$ (and then, of course, to $a$). This gives the pattern
$\gamma'_\rho$.


\define\figtwo{\includegraphics{unim2.ps}}
\midinsert \vbox{\hbox{\hfil\lower3.5in\vbox{\figtwo}\hfill}}
\botcaption {Figure 0.1} Patterns $\gamma_{2/5}$ and $\gamma'_{2/5}$
\endcaption
\endinsert

Also suppose that $\gamma'_{1/2}:\{0,1/4,1/2,3/4\}\to
\{0,1/4,2/4,3/4\}$ is defined as follows: $\gamma'_{1/2}(0)=1/2,
\gamma'_{1/2}(1/4)=3/4, \gamma'_{1/2}(1/2)=1/2,
\gamma'_{1/2}(3/4)=0$. Clearly, the kneading invariant $\nu'_\mu$
corresponding to $\gamma'_\mu$ is obtained from $\nu_\mu$ as
follows: $\nu'_\mu=(\nu_\mu(1), \dots, \nu_\mu(q-2), L, R, R, R,
\dots)$, i.e. in the end of $\nu'_\mu$ there stands an infinite
string of $R$-s.  Let also $\nu'_\mu=\nu_\mu$ if $\mu$ is
irrational.  Finally let us denote the kneading sequence of a
unimodal map $f$ by $\nu(f)$.

\proclaim {Theorem 2.4} Let $f$ be a unimodal map. Then
$I_f=[\mu,1/2]$ iff $\nu'_\mu\succ \nu(f)\succ \nu_\mu$; in particular
if $\mu$ is irrational then $I_f=[\mu,1/2]$ iff $\nu_\mu=\nu(f)$.
\endproclaim

In fact the left endpoint $\rho_f$ of the rotation interval of a
unimodal map $f$ was also introduced as a topological invariant of
the map by J.-M. Gambaudo and C. Tresser in [GT] (see [BM1] where
the connection between the rotation interval and the invariant
introduced in [GT] is established). The authors study the behavior
of $\rho_{f_\nu}$ for families of maps $f_\nu:[-1,1]\to
[-1,1],f_\nu(x)=1-\nu |x|^\mu$ where $\mu\le 1$ and prove for these
families the decreasing of $\rho_{f_\nu}$ which in our terms means
the growth of the rotation interval $I_{f_\nu}$. This fact can be
deduced from the results of Section 3 which deals with some ways to
compare rotation intervals of various interval maps. By Theorem 0.1
the set of periods is defined by the rotation interval if the latter
is not degenerate, so we get a method of comparing both rotation
intervals and sets of periods of maps. As an application we prove
the monotone growth of rotation interval for some one-parameter
families of unimodal maps.

Let $\Cal S=\{f: f$ be a convex map of the interval $[0,1]$ into
itself with a unique turning point $c_f$ which is a local maximum
such that $f|[0,c_f]$ and $f|[c_f,1]$ are $C^1$-maps, $f(c_f)>c_f$
and $f(0)=f(1)=0 \}$.  Let us make some simple remarks. First, we
only require that $f$ be continuously differentiable at $c_f$ from
the left and the right separately. In fact, one-sided derivatives of
$f$ at $c_f$ do not necessarily vanish. Also, the assumption
$f(c_f)>c_f$ is needed to avoid considering trivial cases and can be
made without loss of generality. Under this assumption there is no
fixed point in $(0,c_f]$ and there is a unique fixed point $a_f\in
(c_f,1]$. For any $x\neq c_f$ there is a well-defined point
$x'_f\neq x$ such that $f(x)=f(x'_f)$; also let $c'_f=c_f$.  Also,
$\lambda(K)$ denotes the Lebesgue measure of a measurable set $K$.

\proclaim {Lemma 3.2} The following statements are true. \roster
\item Let $|f'|\ge |g'|, c_f=c_g=c$ and
$\dfrac {|c-a_f|}{|c-a'_f|}\ge \dfrac {|c-a_g|}{|c-a'_g|}$ where
$f,g\in \Cal S$. Then $I_f\supset I_g$.
\item Let $g\in \Cal S$ and $|g'(x)(x-c_g)|\le |g'(x'_g)(x'_g-c_g)|$
for any $x\ge c_g$. Then $I_{\nu g}\supset I_g$ for any $\nu>1$.
\endroster
\endproclaim

Another close result deals with maps from the class $\Cal G\subset
\Cal S, \Cal G= \{g: g\in \Cal S$ is a polynomial map of $[0,1]$
into itself of degree no more than $3 \}$.

\proclaim {Lemma 3.7} Let $f\ge g$ and $f,g\in \Cal G$. Then
$I_f\supset I_g$.
\endproclaim

For the sake of convenience we sum up the results of Lemmas 3.2 and
3.7 in Theorem 3.8 dealing with one-parameter families of interval
maps.

\proclaim { Theorem 3.8} Let $f_\nu, \nu \in [b,d]$ be a
one-parameter family of interval maps such that one of the following
properties holds. \roster
\item $f_\nu\in \Cal S$ for any $\nu$; also,
if $\nu>\mu$ then $|f'_\nu|\ge |f'_\mu|$ and $\dfrac
{|c_{f_\nu}-a_{f_\nu}|}{|c_{f_\nu}-a'_{f_\nu}|}\ge \dfrac
{|c_{f_\mu}-a_{f_\mu}|}{|c_{f_\mu}-a'_{f_\mu}|}$.
\item $f=f_b\in \Cal S,\,|f'(x)(x-c_f)|\le |f'(x'_f)(x'_f-c_f)|$
for any $x\ge c_f$ and $f_\nu=\nu f$.
\item $f_\nu \in \Cal G$ for any $\nu$ and $f_\nu\ge f_\mu$ if $\nu>\mu$.
\endroster

Then $I_{f_\nu}\supset I_{f_\mu}$ if $\nu>\mu$ and so if $f_b$ has
an odd periodic point then $P(f_\nu)\supset P(f_\mu)$ if $\nu>\mu$.
\endproclaim


\head 1. Preliminaries \endhead

We need some well-known tools. Let $I_0,\dots$ be intervals such
that $f(I_j)\supset I_{j+1}$ for $0\le j$; then we say that
$I_0,\dots$ is an {\it $f$-chain} or simply a {\it chain} of
intervals. If a finite chain of intervals $I_0,\dots, I_{k-1}$ is
such that $f(I_{k-1})\supset I_0$ then we call $I_0,\dots, I_{k-1}$
an {\it $f$-loop} or simply a {\it loop} of intervals.

\proclaim{Lemma ALM [ALM2]} The following statements are true.
\roster
\item Let $I_0,\dots,I_k$ be a finite
chain of intervals. Then there is an interval $M_k$ such that
$f^j(M_1)\subset I_j$ for $0\le j\le k-1$ and $f^k(M_k)=I_k$.
\item Let $I_0,\dots$ be an infinite chain of
intervals. Then there is a nested sequence of intervals $M_k$
defined as in {\rm (1)} whose intersection is an interval $M$ such
that $f^j(M)\subset I_j$ for all $j$.
\item Let $I_0,\dots,I_k$ be a loop of intervals. Then
there is a periodic point $x$ such that $f^j(x)\in I_j$ for $0\le
j\le k-1$ and $f^k(x)=x$.
\endroster
\endproclaim

Let $\Cal U$ be the set of all piecewise-monotone interval maps $g$
with one fixed point (denoted $a_g=a$). Fix $f\in \Cal U$; then
$f(x)>x$ for any $x<a$ and $f(x)<x$ for any $x>a$. Call an interval
$I$ {\it admissible} if one of its endpoints is $a$. Call a chain (a
loop) of admissible intervals $I_0,I_1,\dots$ {\it admissible}; if
$I_0,\dots,I_{k-1}$ is an admissible loop then $k>1$ since the image
of an admissible interval cannot contain this interval. For any
admissible loop $\bar \alpha=\{I_0,\dots,I_{k-1}\}$ call the pair of
numbers $(p/2,k)=orp(\bar \alpha)$ the {\it over-rotation pair} of
$\bar \alpha$ where $p$ is the number of indices $0\le s\le k-1$
with $I_s$ and $I_{s+1}$ located on opposite sides of $a$. It is
easy to see, that this definition is consistent with the definition
of over-rotation pair given above. Observe also, that for an
admissible loop the number $p$ is always even (as the interval has
to come back to where the loop starts). Call the number $\rho(\bar
\alpha)=p/2k$ the {\it over-rotation number} of $\bar \alpha$. A
sequence $\{y_1,\dots,y_l\}$ is called {\it non-repetitive} if it
cannot be represented as several repetitions of a smaller sequence.
Define a function $\phi_a$ on all admissible intervals so that
$\phi_a([b,a])=0$ if $b<a$ and $\phi_a([a,d])=1$ if $a<d$. Finally,
given a set $A$ and a point $x$ we say that $A\le x$ if for any
$y\in A$ we have $y\le x$.

\proclaim{Lemma BM3 [BM0]} Let $f\in \Cal U$ and  $\bar
\alpha=\{I_0,\dots,I_{k-1}\}$  be an admissible  loop of
non-degenerate intervals. Then there are the following
possibilities. \roster
\item Let $k$ be even and for each $j$ the intervals $I_j$ and $I_{j+1}$ are
such that either $I_j\le a\le I_{j+1}$ or $I_j\ge a\ge I_{j+1}$.
Then $f$ has a point $x$ of period $2$.
\item If the first possibility fails, then there is a
periodic point $x\in I_0$ such that $x\neq a, f^j(x)\in I_j (0\le
j\le k-1), f^k(x)=x$ and so $\rho(x)=\rho(\bar \alpha)$. If the
sequence $\{\phi_a(I_0), \dots, \phi_a(I_{k-1})\}$ is
non-repetitive, then $orp(x)=orp(\bar \alpha)$. Moreover, $x$ can be
found so that the following holds: for every $y$ from the orbit of
$x$ there exists no $z$ such that $y>_a z$ and $f(y)=f(z)$.
\endroster
\endproclaim

Any point $x$ with the properties from Lemma BM3 is said to be {\it
generated} by $\bar \alpha$.

\proclaim{Lemma 1.1} Let $f\in \Cal U$, let a point $c<a$ be such
that $f(c)>a$ and $f^n(c)\le c$ for some number $n$. Suppose that
among the points $c, f(c), \dots, f^{n-1}(c)$ there are $p$ iterates
$x$ with $x$ and $f(x)$ lying on opposite sides of $a$; then the
following holds.

\roster
\item $p$ is even and $p/2n\in I_f$.
\item If $f^{n+1}(c)\le a$ then for any $r,s$
with $r/s=p/2n$ there is a periodic point of over-rotation pair
$(r,s)$.
\item If $f^{n+1}(c)\le a$ and for some $N$ the $f^N$-image of interval
    $[f^{n+1}(c), a]$ covers $c$ then a small left semi-neighborhood of
    $p/2n$ is contained in $I_f$.
\endroster


\endproclaim

\demo{Proof} (1) As $c$ and $f^n(c)$ are on the same side of $a$, it
follows that $p=2p''$ is even. Consider the admissible loop $[c,a],
[f(c),a], \dots, [f^{n-1}(c),a]$. By Lemma 1.8. it generates a
periodic point $x$ of over-rotation number $p/2n=p''/n$.

(2) Since $f(c)>a$, $f^n(c)\le c$ and $f^{n+1}(c)\le a$ then there
is a point $a'\in [f^n(c), c)$ such that $f(a')=a$.  By Lemma ALM(1)
there is an interval $M\subset [c,a]$ with $fM\subset [f(c),a],
\dots, f^nM=[a', a]$ and an interval $M'\subset [a',c]$ with
$fM\subset [f(c),a], \dots, f^nM'=[a',a]$. Thus, $f^nM\cap
f^nM'\supset M\cup M'$.

Now, let $p''$ and $n$ be coprime. By Lemma
ALM there is a periodic point $x\in M$ with $f^n(x)=x$ and
$\rho(x)=p''/n$. Since $p''$ and $n$ are coprime, then
$orp(x)=(p'',n)$. Using standard arguments, for a given $r>1$ we can
find a periodic point $y\in M$ of period $rn$ with $f^n(y)\in M,
\dots, f^{(r-1)n}(y)\in M, f^{rn}(y)=y$. Then the construction
implies that $orp(y)=(rp'', rn)$. This proves (2) in the case when
$p'', n$ are coprime.

Now, let $p''=tp', n=tn'$ where $p', n'$ are coprime and $t\ge 2$.
By Lemma ALM there is a periodic point $x\in M$ of period $n$;
clearly, it has over-rotation pair $(sp',sn')$ for some $s$. By
Theorem BM2 and by [MN], the pattern $\pi$ represented by $orb(x)$
forces an irreducible pattern $\pi'$ of over-rotation pair
$(p',n')$. Suppose that $\pi$ does not have a block structure over
$\pi'$. Then by [MN] $\pi$ forces the existence of patterns of all
possible over-rotation pairs $(qp', qn'), q\ge 1$ as desired. Thus,
we may assume that $\pi$, the pattern of the orbit of $x$, has block
structure over $\pi'$.

Hence the points from the orbit of $x$ enter intervals $[0,a]$ and
$[a,1]$ periodically with the period $n'$. On the other hand, by the
construction $f^j(c)$ and $f^j(x)\in f^j(M)$ are located on the same
side of $a$ for every $j\le n$. Thus, $f^i(c)$ and $f^{n'}(f^i(c))$
lie to the same side of $a$ for all $0\le i\le n-n'$. Consider
possible locations of $f^{n'}(c)$. By the above $f^{n'}(c)<a$.
Suppose that $f^{n'}(c)<a'$. Then, as in the first paragraph of the
proof of (2), by Lemma ALM(1) we can construct two intervals $N, N'$
(similar to $M, M'$) with $N\subset [c,a], fN\subset [f(c),a],
\dots, f^nN=[a', a]$ and $N'\subset [a',c], fN'\subset [f(c),a],
\dots, f^nN'=[a',a]$. This implies that for any $r>1$ we can find a
periodic point $y\in N$ with $f^{n'}(y)\in N, \dots,
f^{(r-1)n'}(y)\in N, f^{rn'}(y)=y$ so that $orp(y)=(rp', rn')$ as
desired. Thus we may assume that $f^{n'}(c)>a'$.

Choose the greatest image $d=f^{n'j}(c)$ of $c$ under powers of
$f^{n'}$ such that $a'<d$. Then $n'j\le n-n'$ (because $f^n(c)<a'$
by the construction) and $f^{n'}(d)\le a'<d$. By the periodicity
with which the orbit of $c$ enters intervals $[0,a]$ and $[a,1]$ we
see that among points $d, f(d), \dots, f^{n'}(d)=f^n(c)\le a'$ there
are $p'$ points lying to the right of $a$. By repeating the
construction from the first paragraph of the proof of (2) and using
Lemma ALM(1) one can find intervals $N, N'$ with $N\subset [d,a],
fN\subset [f(d),a], \dots, f^nN=[a', a]$ and $N'\subset [a',d],
fN'\subset [f(d),a], \dots, f^nN'=[a',a]$. So repeating the
arguments from the first paragraph of the proof we can find points
of all over-rotation pairs $(rp',rn')$.

(3) Let $b\neq a$. If $b<a$ then $f_r(b)=\max \{f(z): z\in [b,a]\}$;
if $b>a$ then $f_r(b)=\min \{f(z): z\in [a,b]\}$. Clearly $f_r$ maps
points from $[0,a]$ into $[a,1]$ and vice versa. Then $[a,
f_r(b)]=f[b,a]\cap [a,1]$ if $b<a$ and $[f_r(b), a]=f[a,b]$ if
$a<b$. Also, since $f\in \Cal U$ then $f([b, a])\subset [b, f_r(b)]$
if $b<a$ and $f([a, b])\subset [f_r(b), b]$ if $a<b$. In the
situation of the lemma there is the smallest $N$ such that
$f_r^N(f^{n+1}(c))\le c$. Moreover, by the properties of $f_r$ we
see that $N=2m$ is even. Hence by the definition of $f_r$ we see
that the following is an admissible loop:
$$[c, a], [f(c),a], \dots, [f^{n+1}(c), a], [f_r(f^{n+1}(c)), a],
\dots, [f_r^{2m-1}((f^{n+1}(c)), a].$$ A direct computation shows
that its over-rotation number is $\frac{p/2+m}{n+2m}<p/2n$. Since by
Lemma BM3 there exists a periodic point with the over-rotation
number $\frac{p/2+m}{n+2m}$, then the proof is complete.
{\unskip\nobreak\quad$\square$} \enddemo

Finally we state the results of [BM3]. One of them gives a criterion
for a pattern to be an over-twist pattern. To state this criterion
we need to define a {\it code}, i.e. a special function which maps
points of either a periodic orbit or of a pattern to the reals. We
also need a few other definitions.

Recall, that a cycle (and the pattern it represents) is {\it divergent\/} if it
has points $x<y$ such that $f(x)<x$ and $f(y)>y$, and that a cycle (pattern)
that is not divergent is called {\it convergent\/}. Clearly, a
convergent pattern has a unique complementary interval $U$
such that its left endpoint is mapped to the right and its right
endpoint is mapped to the left. If this pattern is exhibited by a
map $f$, then this interval contains a fixed point, always denoted
by $a$. Given two points $x, y$ of the pattern we say that $x>_a y$
if $x$ and $y$ are located on the same side of $U$ and $x$ is
farther away from $U$ than $y$.

Similar notation is used for periodic orbits of interval maps which
exhibit convergent patterns. There is an equivalent way to define
convergent patterns. 
Namely, if $f$ is a $P$-monotone map for a cycle $P$
then $P$ is convergent if and only if $f\in\Cal U$. If $f\in \Cal U$
and $a$ is the fixed point of $f$ then we write $x>_a y$ if points
$x, y$ are located on the same side of $a$ and $x$ is farther away
from $a$ than $y$. (This notation is similar to the one used for
convergent patterns.)

Let $P$ be a cycle of $f\in\Cal U$ and $\varphi$ be a function
defined as $1$ to the right of $a$ and zero elsewhere. Following
[BK], we introduce the {\it code\/} for $P$ as follows. The code is
a function $L:P\to\R$, defined by $L(x)=0$ for the leftmost point
$x$ of $P$ and then by induction we have
$L(f(y))=L(y)+\rho-\varphi(y)$, where $\rho$ is the over-rotation
number of $P$. When we get back to $x$ along the orbit $P$, we add
$\rho$\ $n$ times ($n$ is the period of $P$), and we subtract the
sum of $\varphi$ along $P$, which is $n\rho$, so we have a sum of
$0$. Therefore, the definition is correct.

Clearly, we can also speak of codes for patterns. If $f\in\Cal U$
(or if the pattern in question is convergent), we say that the code
for $P$ is {\it monotone\/} if for any $x,y\in P$, $x>_a y$ implies
$L(x)<L(y)$.

\proclaim{ Theorem BM4 [BM3]} A pattern is over-twist if and only if
it is convergent and has monotone code.
\endproclaim

\head 2. Unimodal over-twist patterns \endhead

Let us describe unimodal over-twist patterns. Our aim is to show
that the pattern $\gamma_\rho$ defined in Subsection 0.3 is the
unique unimodal over-twist pattern of over-rotation number $\rho$.

\proclaim{ Lemma 2.1} Suppose that $\rho=p/q$ is such that
$0<\rho\le 1/2$. Then the only unimodal over-twist pattern of
over-rotation number $\rho$ is the pattern $\gamma_\rho$. Thus, a
unimodal map $f$ has the over-rotation interval  $[\mu, 1/2]$ with
$\mu\le \rho$ if and only if $f$ has a periodic orbit of pattern
$\gamma_\rho$.
\endproclaim

\demo{Proof} Suppose that $\tau$ is a unimodal over-twist pattern of
over-rotation number $\rho=p/q<1/2$. Consider a unimodal map $f:[0,
1]\to [0, 1]$ with the unique critical point $c$, unique fixed point
$a$, and the unique $f$-preimage $a'$ of $a$. Assume that $c$ is
periodic and its orbit $P$ exhibits the pattern $\tau$. Then by
Lemma BM4 there are no points of $P$ in $[a', c)$. Moreover, by
Theorem BM2 $c$ is of period $q$. By definition of over-rotation
pair, there are $p$ points of $P$ in the interval $(a, 1]$. This
implies that there are $p$ points of $P$ in the interval $[c, a]$
too and, moreover, that $f(P\cap [c, a])=P\cap [a, 1]$. In fact, the
map $f$ simply flips all points of $P\cap [c, a]$ to the other side
of $a$. Clearly, then there are $q-2p$ points of $P$ in $[0, c)$.

Let us compute out codes of some points of $P$. By definition
$L(0)=0$. This implies that $L(1)=1-\rho$ and
$L(c)=1-2\rho=(q-2p)/q$. Since the values of the code on points of
$P$ are fractions with denominator $p$ and the code on $P\cap [0,
c)=\{x_0=0<x_1<\dots<x_{q-2p-1}\}$ is monotonically increasing, we
see that on all the $q-2p$ points of $P\cap [0, c)=P\cap [0, a']$
the code equals $L(x_0)=0, L(x_1)=1/q, \dots,
L(x_{q-2p-1})=(q-2p-1)/q$. This (and the definition of the code)
immediately implies that $\tau=\gamma_\rho$ as desired.

Now, by Theorem Bo1, Theorem B1 and Theorem BM2 it follows that if a
unimodal map $f$ has the over-rotation interval $I_f=[\mu, 1/2]$ and
$\mu\le \rho$ then $f$ has a periodic orbit of pattern
$\gamma_\rho$. On the other hand, it is clear that if $f$ has a
periodic orbit of pattern $\gamma_\rho$ then $I_f=[\mu, 1/2]$ with
$\mu\le \rho$. This completes the proof of the lemma.
{\unskip\nobreak\quad$\square$} \enddemo

Lemma 2.1 gives the least (the weakest) kneading invariant of $f$
implying the fact that $\rho\in I_f$. It also follows from Lemma 2.1
and Theorem BM2 that the kneading invariant $\nu_{p/q}$ of the
over-twist pattern of rational rotation number $p/q$ depends
monotonically on $p/q$: if $s<t$ then $\nu_s\succ \nu_t$. This fact
can be also easily checked directly. However, using kneading
sequences is much more to the point if we concentrate upon
the question of the greatest (strongest) kneading sequence (pattern)
which gives the over-rotation interval precisely coinciding with
 $[\rho, 1/2]$.

\proclaim {Lemma 2.2} A unimodal kneading sequence $\nu$ is such
that $I_\nu=[p/q, 1/2]$ if and only if $\nu'_{p/q}\succeq \nu
\succeq \nu_{p/q}$.
\endproclaim

\demo{Proof} By Lemma 2.1 if $\nu_{p/q}\le \nu$ then $[p/q]\subset
I_\nu$. Let us show that $I_{\nu'_{p/q}}=[p/q, 1/2]$. Indeed, assume
that $f$ is a unimodal map such that its unique critical point $c$
eventually maps to its unique fixed point $a$ so that the orbit of
$c$ exhibits the pattern $\gamma'_{p/q}$. We need to show that
$I_f=[p/q, 1/2]$. By Lemma 2.1 it suffices to prove that for any
periodic orbit $Q$ of $f$ we have that $\rho(Q)\ge p/q$.

There are two ways this claim can be proven. One of them uses
kneading sequences. Indeed, suppose that for some cycle $Q$ of $f$
we have that $\rho(Q)<p/q$. By Lemma 2.1 we may assume that $Q$ has
an over-twist pattern. For convenience we can choose a very big $r$
and then choose $Q$ of period $rq+1$ and over-rotation number
$rp/(rq+1)$. Then it follows from the description of the pattern
$\gamma'_{p/q}$ and from the properties of $Q$ given by Lemma 2.1,
that $i(x)\succ \nu(f)$, a contradiction with
properties of kneading sequences.

However, we will give a more direct proof of the claim. Assume that
$P$ is a periodic orbit of pattern $\gamma_{p/q}$ with $\min P=0$
and $\max P=1$. Let $f$ be the $P$-linear map. Construct the
oriented graph $G_P$ whose vertices correspond to the closures of
the components of $[0, 1]\setminus \{P\cup a_f\cup a_f'\}$ and whose
arrows connect a vertex $J$ and a vertex $I$ if and only if
$f(J)\supset I$ (this is a standard construction in one-dimensional
dynamics [ALM2]; as always, $a_f$ is the fixed point of $f$ and
$a_f'$ is its $f$-preimage). Assign numerical value $1$ to all
arrows which come out of segments-vertices located to the right of
$a$, and assign $0$ to all arrows which come out of segments-vertices
located to the left of $a$. Denote the just defined function on
arrows by $\psi$. Then it follows that the over-rotation number of
$P$ is the minimum of averages of assigned values of $\psi$ along
all possible loops in the graph $G_P$.

Now, suppose that $Q$ is the orbit of (non-cyclic) pattern
$\gamma'_{p/q}$ with $\min Q=0$ and $\max Q=1$; denote the map
acting on $Q$ by $g$. We may assume that points of $Q$ coincide with
points of $P$ except for the fact that the point $z$ of $P$ which
maps to $c_f$ by $f$ will be mapped to $a'_f$ by $g$ (the map $g$
will act the same way on $a'_f, a_f$ as $f$). This allows us to use
notation $a, a', c$ without references to the map $f$. Extend $g$ to
the $Q$-linear map. Then construct the oriented graph whose vertices
correspond to the closures of the components of $[0, 1]\setminus
\{Q\cup a_g\cup a_g'\}$ and whose arrows connect a vertex $J$ and a
vertex $I$ if and only if $g(J)\supset I$ ($a_g$ is the fixed point
of $g$ and $a_g'$ is its $g$-preimage). Clearly, the sets of
vertices of the graphs $G_P$ and $G_Q$ are the same.

Assign numerical value $1$ to all arrows of $G_Q$ which come out of
segments-vertices located to the right of $a_g$ and $0$ to all
arrows of $G_Q$ which come out of segments-vertices located to the
left of $a_g$. It follows that the graphs $G_P$ and $G_Q$ almost
coincide, except for arrows of either graph which come into the
segment-vertex $[a', c]$. Now, suppose that $g$ has a periodic orbit
$Z$ of over-rotation number less than $p/q$. By Lemma BM3 we may
assume that $Z$ avoids $[a, c]$. However this implies that the loop
of arrows in $G_Q$ which corresponds to $Z$ avoids arrows which come
into the segment-vertex $[a', c]$. Hence this loop of arrows
consists of the arrows common for both $G_P$ and $G_Q$. This is a
contradiction as all loops of arrows in $G_P$ must produce averages
of $\psi$ which are greater than or equal to $p/q$.

It remains to show that if a kneading sequence $\nu$ is such that
$\nu\succ \nu'_{p/q}$ then $I_\nu=[t, 1/2]$ and $t<p/q$. Suppose
that $\nu\succ \nu'_{p/q}$ and denote a unimodal map which exhibits
this itinerary by $h$. We can find a point $d$, say, to the right of
$c_h$ so that except for the first moment the $h$-itinerary of $d$
and the $g$-itinerary of $c_g$ under the above defined map $g$ are
the same up to the point when $c$ maps to $a'$ by $g^q$ and $d$ maps
slightly to the left of $a'$ by $h^q$; now Lemma 1.1(3) immediately
implies that $I_h$ contains not only $p/q$ but also its small
neighborhood, a contradiction. {\unskip\nobreak\quad$\square$}
\enddemo

Lemma 2.2 gives the precise description of unimodal maps which have
the over-rotation interval $[p/q, 1/2]$ through their kneading
sequences. However, we still do not have the description of unimodal
maps which have the over-rotation intervals with {\it irrational}
left endpoint. We will describe them in Lemma 2.3. In the proof we
use a special construction which relates unimodal maps to irrational
rotations of the circle.

\proclaim {Lemma 2.3} Let $f$ be a unimodal map. If $\mu$ is
irrational then $I_f=[\mu,1/2]$ iff $\nu_\mu=\nu(f)$.
\endproclaim

\demo{Proof} Consider a unimodal map $f$ such that $I_f=[\rho_f,
1/2]$ and $\rho_f$ is irrational. Then, clearly, the kneading
sequence of $f$ is not periodic or preperiodic. It is well-known
that then we may assume that $f$ has no wandering intervals. If $f$
had such wandering intervals, they could be collapsed to points which would not
change the over-rotation interval. Similarly, we may assume that
$f$ has no non-trivial periodic intervals. Otherwise either the intervals could
be collapsed to points, or they contain $c_f$ which implies that
$\rho_f$ is rational. This in turn implies that $c_f$ is approached
from either side by periodic or preperiodic points.

Choose any periodic or preperiodic point $z$ and then choose the
point $x$ in the orbit of $z$ so that $f(x)$ is greater than any
point of the orbit of $z$. Consider a {\it truncation} $f_x$ of $f$
defined as follows: for every $y$ we have that $f_x(y)=f(y)$ if
$f(y)\le f(x)$ and $f_x(y)=f(x)$ if $f(y)>f(x)$. Then $f_x$ has a
preperiodic critical point which implies that it has the
over-rotation interval with rational left endpoint. Hence this
interval is strictly smaller than $I_f$. We conclude that there are
kneading sequences which are smaller than $\nu_f$ and generate
smaller over-rotation intervals that are arbitrarily close to $I_f$.
A similar analysis shows that the change in $I_f$ can be achieved
also if we increase $\nu_f$. In this respect, the over-rotation
interval depends upon kneading sequence sensitively at $\nu_f$.

In fact, for a given irrational number $\rho<1/2$ there is a unique
kneading sequence $\nu_\rho$ such that $I(\nu_\rho)=[\rho, 1/2]$.
This kneading sequence is defined earlier in Subsection 0.3 and is
closely related to circle rotations (i.e., to shifts by $\rho \mod
1$ considered on $[0, 1)$. Thus, the kneading sequences $\nu_\rho$
can be characterized as the $\succ$-smallest such that the
corresponding over-rotation interval $I(\nu_\rho)$ contains the
number $\rho$.

Let us verify the claims made in the previous paragraph. Unlike in
the case of rational numbers $\rho$ where we relied upon interval
techniques, for irrational numbers $\rho$ we develop a new approach
directly relating unimodal maps $f$ to circle maps. Here we use a
special discontinuous lifting of $f$ to a discontinuous degree one
map of the real line, studied in the spirit of [M3]. The
construction is described below. One can easily see that using this
construction we can treat both rational and irrational cases.
However we chose to use interval tools to tackle the rational case
to show how different methods can work with unimodal maps.

Let us fix a unimodal map $f$ with the rotation interval
$[\mu,1/2], \mu<1/2$; it implies that $f(0)<a$ because otherwise we
have that $f([0, a])\subset [a, 1], f([a, 1])\subset [0, a]$ and
thus $I_f=\{1/2\}$, a contradiction. As we have already considered
the case when $\mu$ is rational, we may assume that $\mu$ is
irrational. Then the construction is as follows.

First we define a discontinuous conjugacy $\sigma:[0,1]\to
[0,1]$ as the identity on $[0,a)$ and the symmetry (flip) with respect to
$(a+1)/2$ on $[a,1]$, so that $\sigma(x)=x$ if $0\le x<a$ and
$\sigma(x)=a+1-x$ if $a\le x\le 1$; note that $\sigma^{-1}=\sigma$.
Define now $g:[0,1]\to [0,1]$ by $g=\sigma\circ f\circ\sigma$ with
the following two changes: $g(a')=a$, not $1$, and $g(1)=a$, not $1$
(see Figure~2.1). Note that $g$-image of $[0,1]$ is $[0,1)$, i.e. it
does not contain $1$.

    1) On the interval $[0,a']$ we have $g(x)=f(x)$, thus the
graph of $g$ on $[0,a']$ is the same as that of $f$; in particular
$g(0)=f(0)$ and $g(a')=a$.

    2) On the interval $(a',a)$ we have
$g(x)=\sigma(f(x))=a+1-f(x)$. We obtain the graph of $g$ by flipping the
corresponding piece of $f$ in the vertical direction symmetrically with
respect to the horizontal line $y=(a+1)/2$; in particular,
$g(c)=a,$ and $g$-images of points that are close to $a$ from the left will approach $1$
and $g$-images of points close to $a'$ from the right will approach $1$.

    3) On the interval $[a,1]$ we have
$g(x)=f(\sigma(x))=f(a+1-x)$, thus the graph of $g$ on $[a,1]$ is
obtained by flipping the corresponding piece of the graph of $f$ in
the horizontal direction symmetrically with respect to the line
$x=(a+1)/2$; in particular $g(a)=0, g(1)=a$.

    Clearly we can define over-rotation numbers for $g$ like it is
done for $f$ (e.g., counting how many times the $g$-orbit of a point
enters the interval $[a,1]$); moreover, the same way the over-rotation
pairs for $g$-periodic orbits may be introduced. Moreover, $\sigma$
conjugates $f$ and $g$ on orbits which avoid $a$ and $c$ so that the
over-rotation sets on these orbits coincide; the same can be said
about over-rotation pairs of periodic points.  Consider the remaining
orbits. First assume that $c$ is neither periodic nor mapped into
$a$ by some iterate of $f$; it means that $g$-orbit of $0$ never
passes through $a$ or $1$ and so $\sigma$ conjugates $f$ on the
$f$-orbit of $0$ with $g$ on the $g$-orbit of $0$ keeping the
rotation sets. Since $g(c)=a, g(a)=0$ we now see that in fact due to
the conjugacy $\sigma$ the rotation sets of $f$ and $g$ are the
same. Similarly we can easily check that the over-rotation numbers
of points for $g$ and of $f$ coincide in the remaining cases (when
$c$ is periodic or $c$ is a preimage of $a$) too. We conclude that
$I_g=I_f=[\mu,1/2]$.

    Now let us specify a lifting $F$ of a map $f$.  We do it by
setting on $[0,1)$
$$F(x)=\cases {g(x)=f(x)}&\text{if $0\le x\le a'$,}\\
{g(x)=a+1-f(x)}&\text{if $a'<x<a$,}\\ {g(x)+1=f(a+1-x)+1}&\text{if
$a\le x< 1$,}\endcases$$ and then as usual: if $x=k+y$ with
$y\in[0,1)$ then $F(x)=k+F(y)$. The map $F$ is a {\it degree one
lifting} map of the real line into itself or an {\it old} map (see
\cite{M3}). Obviously by the construction the sets of classical
rotation numbers and pairs of $F$ coincide with the sets of
over-rotation numbers and pairs of $g$, and
therefore with the sets of
over-rotation numbers and pairs of $f$, so the classical rotation set $I_F$
coincides with $I_f$. An example can be found on Figure 2.1.

\define\figthree{\includegraphics{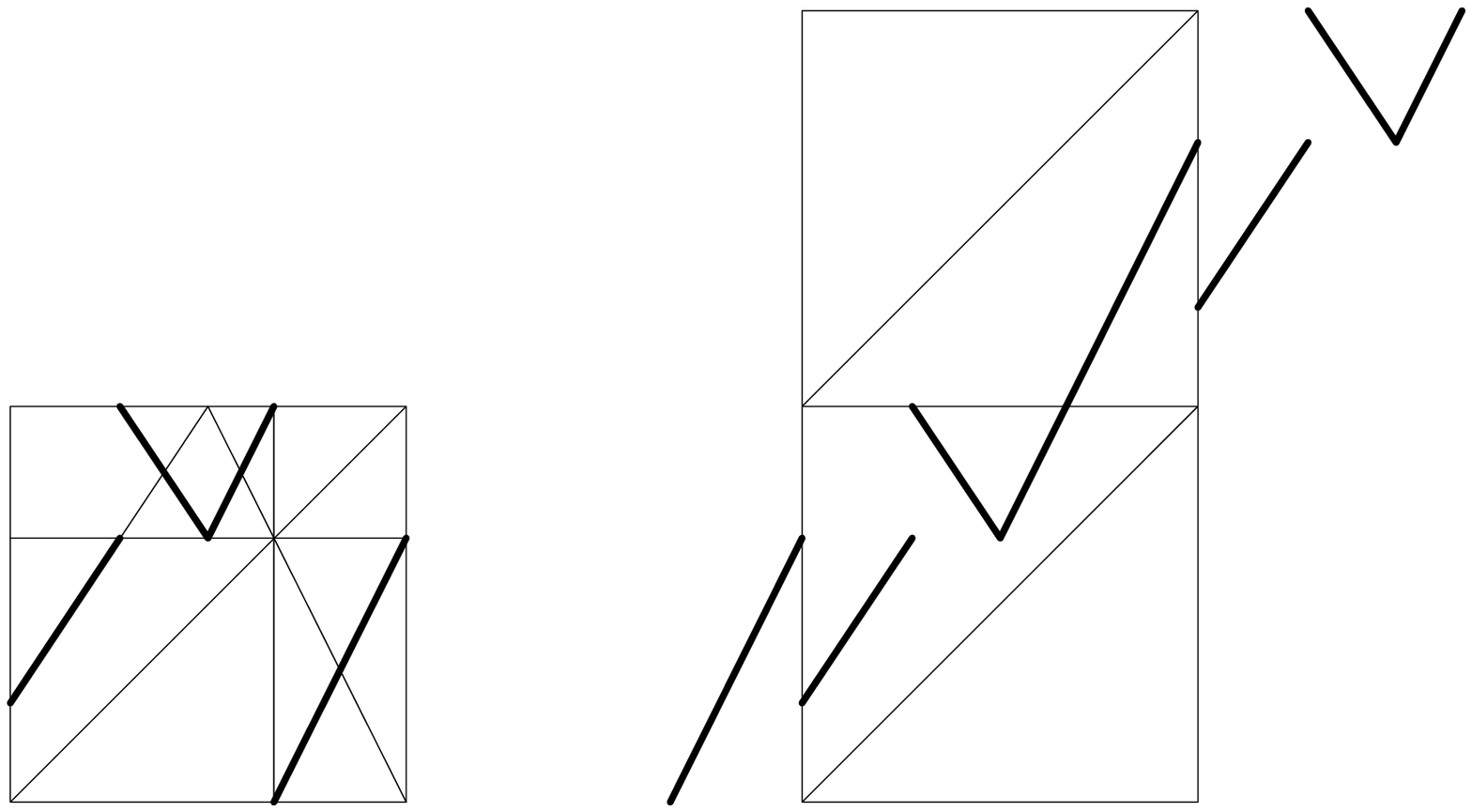}}
\midinsert \moveright 0.375in \vbox{\hbox
{\hfil\lower5.2in\vbox{\figthree}\hfil}} \vskip-.4in \botcaption
{Figure 2.1} Functions $g$ (left) and $F$ (right)
\endcaption \endinsert

A famous nice tool for studying of old maps is so-called ``pouring
water from below or above'' (see [M1-M3]); the definitions can be given in
general situation but we only explain how one can use ``pouring
water'' in our particular case. Let $d\ne 0$ be such that
$f(0)=f(d)$).  Consider the continuous map $G:\Bbb R\to \Bbb R$ of
degree one which depends on $F$ (i.e. in the end on $f$) and is
defined on every interval $[n,n+1]$ as follows: $$G(x)=\cases
{F(x)}&\text{if
$n\le x\le n+a'$,}\\ {n+a}&\text{if $n+a'<x\le n+c$}\\
{F(x)}&\text{if $n+c<x\le n+a+1-d$}\\ {F(n+1)=n+f(0)}&\text{if
$n+a+1-d<x<n+1$}\endcases$$ The connection between $F|[n,n+1]$ and
$G|[n,n+1]$ is obvious: (1) $F=G$ except for two intervals, $(n+a',
c)$ and $(n+a+1-d, n+1)$ on which $G$ is a constant; (2) $G\le F$;
(3) $G$ is continuous.  Clearly the construction is possible for a
unimodal map with homtervals; let us see though what consequences
the absence of homtervals implies for $G$. First, it follows that
the only flat spots of $G$ are intervals $(n+a', c)$ and $(n+a+1-d,
n+1)$. Now, suppose there is an interval $I\in [0,1]$ such that
$G^m(I)$ is disjoint from intervals $(n+a', c)$ and $(n+a+1-d, n+1)$
for any $m$ and $n$. The connection between $f$ and $G$ then implies
that $f^m|I$ is monotone for any $m$ which is impossible since $f$
has no homtervals.  So for any interval $I$ there are $m$ and $n$
such that $G^m(I)$ is not disjoint from $(n+a', c)\cup (n+a+1-d,
n+1)$. In other words the set $A$ of all points $X$ which avoid all
intervals $(n+a', c)$ and $(n+a+1-d, n+1)$ is nowhere dense.

    The properties of continuous monotone old maps of real line into
itself are studied in [M1-M3], [ALM2] (for the sake of convenience
in what follows we mostly refer to [ALM2] with respect to related
questions). It is proven there that there is a unique number $\rho$
such that for any $Z\in \Bbb R$ we have $(1/n) G^n(Z)\to \rho$ and
that there is a point $Z'\in [0,1)$ whose $G$-orbit avoids intervals
$(n+a', c)$ and $(n+a+1-d, n+1)$ for all $n$. Therefore due to the
connection between $f$ and $G$ we have $I_G(Z')=\rho=I_f(Z')$ which
implies that $\mu\le \rho$. On the other hand $\rho\le \mu$ since
$G\le F$ and $G$ is monotone; indeed, for any point $X\in \Bbb R$ we
have that $G(X)\le F(X)$ and moreover, $G^n(X)\le F^n(X)$ implies
(after we apply $G$ to both sides of the inequality) that
$G^{n+1}(X)\le G(F^n(X))\le F^{n+1}(X)$. Hence when we take the
limits of $G^n(X)/n$ and $F^n(X)/n$ we see that $\rho$ is less than
or equal to all numbers from the over-rotation interval of $f$, i.e.
that $\rho\le \mu$ (see, e.g., [ALM2]). So, $\rho=\mu$.

Let $\pi:\Bbb R\to S^1$ be the usual projection of $\Bbb R$ onto
$S^1$ such that $[0,1)$ is mapped onto the circle $1$-to-$1$.  Then
$\pi$ semiconjugates $G$ to a continuous monotone map $\ftil:S^1\to
S^1$ of degree $1$ of the same rotation number $\mu$. Moreover,
$\ftil$ has exactly two arcs which it collapses to points and
otherwise is strictly monotone. These arcs are $J=\pi[a+1-d,1]$ and
$V=\pi[a',c]$. Clearly, $\pi(Z')$ is a point of the circle which
avoids interiors of these arcs. The closure $H$ of the orbit of
$\pi(Z')$ in the circle is then also disjoint from the interiors of
$J$ and $V$. It is well-known, that the induced map on $H$ can be
semiconjugate to the interval rotation by the angle $\mu$ by means
of collapsing complementary to $H$ arcs. Moreover, this implies that
$H$ is minimal.

Let us show that $H$ must contain $\pi(c)$ (so that $c\in
\pi^{-1}(H)$). Indeed, suppose otherwise. Then it follows from the
construction that there is an invariant closed set $H'\subset [0,
1]$ which corresponds, through the construction, to the set $H$ in
the circle. Using the above introduced notation we see that $H'$
consists of points which do not enter $(a', c)$ and $(a, d)$. In
other words, points of $H'$ avoid $(a', c)$ and those of them which
belong to $(a, 1)$ have images located non-strictly to the left of
$f(0)$. By construction, points $x'$ of $H'$ have the same
over-rotation numbers as the limits of $G^n(x)/n$ for the points
$x\in H$, i.e. $\rho$.

If, by way of contradiction, $\pi(c)$ does not belong to $H$, then
it follows that $c$ does not belong to $H'$. Hence, as was explained
right after the proof of Lemma 2.2 we can find periodic and
preperiodic points $x$ as close to $c$ as we wish so that these
points will never have eventual images greater than $f(x)$. Given
such point $x$ we can construct a well-defined truncation $f_x$ of
$f$ so that all the points which never enter the open segment
$(f(x), f(c))$ have the same orbits under both $f$ and $f_x$.
Moreover, we can choose $x$ arbitrarily close to $c$. In particular,
we can do this so that $f|_{H'}=f_x|_{H'}$ . Then on the one hand
the over-rotation interval of $f_x$ must be such that
$\rho_{f_x}>\mu=\rho$ is rational and hence $\rho_{f_x}>\mu$, on the
other hand points of $H'$ produce over-rotation numbers
$\rho=\mu<\rho_{f_x}$, a contradiction. Hence $c\in H'$ as desired.
Since $H$ is minimal, it follows that so is $H'$. In particular,
$H'=\omega(c)$. Again by construction this implies that the kneading
sequence of $f$ coincides with the kneading sequence $\nu_\mu$
defined in Subsection 0.3.

To complete the proof of the lemma, it remains to show that for any
irrational number $\mu$ there exists a unimodal map $f$ which has
the over-rotation interval $[\mu, 1/2]$. To construct such a map we
reverse the construction. First, we construct a monotone map $\ftil$
of the circle to itself which has two ``flat spots'', i.e. two arcs
$V=[a', c]$ and $J=[a+1-d, 1]$ which $\ftil$ collapses to points and
$0<a'<c<a<a+1-d<1$ are such points of the circle that $\ftil(a')=a$
and $\ftil(a)=1$. Moreover, the point $0=1$ of the circle never
enters $V\cup J$ and is such that the order of points in its
$\ftil$-orbit is the same as the order of points in the orbit of a
point of the circle under the irrational rotation by the angle $\mu$
(it is easy to see that this is possible). This implies that the
rotation interval of $\ftil$ (as defined in [ALM2]) is degenerate
and coincides with $\mu$. It remains now to reverse the construction
in order to see that the corresponding to $\ftil$ unimodal map $f$
exists. Reversing the arguments from above we also see that
$I_f=[\mu, 1/2]$ as desired. By definition and by construction the
kneading sequence of $f$ coincides with the kneading sequence
$\nu_\mu$ defined in Subsection 0.3. This completes the proof of the
lemma. {\unskip\nobreak\quad$\square$}
\enddemo

It is now easy to see that if we put together Lemmas 2.1, 2.2.and
2.3 we get Theorem 2.4.

\proclaim {Theorem 2.4} Let $f$ be a unimodal map. Then
$I_f=[\mu,1/2]$ iff $\nu'_\mu\ge \nu(f)\ge \nu_\mu$; in particular
if $\mu$ is irrational then $I_f=[\mu,1/2]$ iff $\nu_\mu=\nu(f)$.
\endproclaim

The construction from Lemma 2.3 shows that if $f$ is a unimodal map
such that $I_f=[\mu, 1/2]$ and $\mu$ is irrational then
$f|_{\omega(c)}$ is semiconjugate to the irrational rotation by the
angle $\mu$ by a map constructed in the proof of Lemma 2.3. More
precisely, to construct this semiconjugacy we need to apply the map
$\sigma$ from the proof of Lemma 2.3 to $\omega(c)$, then transport
$\sigma(\omega(c))$ to circle using the standard projection $\pi$ of
the interval onto the circle, and then consider the map induced by
$f$ on the set $\pi(\sigma(\omega(c)))$. The resulting map of a
closed subset of the circle will be the map $\ftil$ from the proof
of Lemma 2.3 and will be such that the order of points in the orbit
of any point of this set is the same as the order of points under the
(irrational) rotation by the angle $\mu$.

On the other hand, it can be verified directly that the same holds
if the orbit of $c$ is periodic and exhibits an over-twist pattern.
That is, in this case we can construct the map $\ftil$ as well;
clearly, it will be defined on a finite set of points, cyclically
permuted by $\ftil$. Lemma 2.1 would imply that then the map $\ftil$
permutes the points of this set as the (rational) rotation by $\mu$
prescribes. Thus, the construction from Lemma 2.3 yields the
description of over-twist patterns too.

In fact, the same construction shows how we can figure out the
number $\rho_f$ for a given unimodal map $f$. Apply the construction
from Lemma 2.3 and construct the map $\ftil:S^1\to S^1$; moreover,
let us use the notation from the proof of Lemma 2.3. Then again, as
in the proof of Lemma 2.3, by [ALM2] there exists a point
$Z'=\pi(z')$ which avoids $V$ and $J$ so that $z'$ avoids $[a', c]$
and $[a, d]$. By [ALM2] the limit $\rho$ of $G^n(Z')/n$ equals the
limit of $G^n(x)/n$ for all $x$; this common limit $\rho$ is called
the {\it rotation number} of $\ftil$. By construction it follows
that the over-rotation number of $z'$ equals $\rho$. Similar to the
proof of Lemma 2.3 one can show that then $I_f=[\rho, 1/2]$ which
implies that $\rho=\rho_f$.

Hence the algorithm of finding $\rho_f$ is as follows. Take the
following three intervals: $K_1=[0, a'], K_2=[c, a]$ and $K_3=[d,
1]$ where $d$ is chosen so that $f(0)=f(d)$. These intervals have
images which complement each other to the whole $[0, 1]$ and
intersect only at the images of their endpoints. Moreover, on each
interval the map is monotone. This implies that there exist points
$x$ with orbits contained in the union $K$ of intervals $K_1, K_2$
and $K_3$. In fact one such point is a fixed point $a$. Moreover, we
can also take the point $a'$ or, more generally, other points which
travel within $K=K_1\cup K_2\cup K_3$ and eventually map to $a$.

However we need other points which travel inside $K$. To discard $a$ and
its preimages, we choose a point $d'\in [c, a]$ such that
$f(d')=d$. Then set $K'_2=[c, d']$, and consider all points $x$ which
travel inside $K'=K_1\cup K'_2\cup K_3$. Since $K'\subset f(K')$,
such points exist. Moreover, by construction it follows that each
such point can play the role of $z'$ exhibits the minimal over-rotation
number $\rho_f$ in $I_f$. This
gives a useful algorithm of finding $\rho_f$ as well as orbits on
which $\rho_f$ is realized as the over-rotation number.

\head 3. How to compare rotation intervals of maps? \endhead

    It is sometimes important to be able not only to estimate or
compute various characteristics of maps but also to compare them for
different maps without in fact estimating or computing. It turns out
that the problem of comparing over-rotation intervals of two
continuous interval maps (not necessarily piecewise monotone) can be
approached from the point of view related to the tools introduced in
[B1] which involve such notions as chains and loops of admissible
intervals (see Lemmas ALM, BM3 and 1.1). In fact a very simple
geometrical condition allows to compare the over-rotation intervals of
two maps.

Suppose that $I\subset J$ and $f:J\to J$ and $g:I\to I$ are two
maps. Let $g$ have a fixed point $a$ and $[f(x),a]\supset [g(x),a]$
for any $x\in I$. Then we say that $f$ is more {\it repellent from
$a$} than $g$. Then if $g(x)<a$ then $f(x)<a$, and if $g(x)>a$ then
$f(x)>a$. Hence if $g(y)=a$ and $y$ is not a local extremum of $g$,
then arbitrarily close to $y$ there are points mapped by $g$ (and
therefore, $f$) both to the left and to the right of $a$. This
implies that $f(y)=a$ too. In particular, if we assume that $a$
itself is not a local extremum of $g$, 
then $f(a)=a$.

    There is also another remark worth making. So far in the
definition we compare how maps $f$ and $g$ repel {\it all} the
points of the interval $I$ from $a$. It is sometimes useful to
compare the repelling of the points of some specific set $B\subset
I$ only. So if $B\subset I$ and we know that $[f(x),a]\supset
[g(x),a]$ for any $x\in B$ then we say that $f$ is more {\it
repellent from $a$ on $B$} than $g$. In fact we deal with the
following specific set $B$. Let $g$ be a map with a unique fixed
point $a$ which is not a local extremum of $g$.
For any point $x$ let $\psi_g(x)$ be the closest to $a$
point in $[x,a]$ (or in $[a,x]$ if $a<x$) such that
$g(x)=g(\psi_g(x))$.
The set $B$ on which we compare the repelling
properties of $f$ and $g$ is the set $\psi_g(I)$. In the statement
of Lemma 3.1 we use the notation introduced above.

\proclaim {Lemma 3.1} The following statements are true.

\roster

\item Let $f:J\to J$ be more repellent from $a$
than $g:I\to I\subset J$ on the set $B=\psi_g(I)$ ($a$ is a unique
fixed point of $g$ inside $I$ which is not a local extremum of $g$). 
Then $I_f\supset I_g$; in particular if $g$ has a
point of odd period then $P(f)\supset P(g)$.

\item Let $g:[0,1]\to [0,1]$ and $f:[u,v]\to [u,v]\subset [0,1]$ be
two maps. Suppose that $g(0)=g(1)=0, g|[0,c]$ is increasing,
$g|[c,1]$ is decreasing, $a\in (c,1]$ is a unique $g$-fixed point in
$(0,1]$, $a'$ is a unique $g$-preimage of $a$ in $[0,c]$,
$[g^2(c),g(c)]\subset [u,v]$ and $f\le g$ on $[g^2(c),a']$, $f\ge g$
on $[c,a]$, $f\le g$ on $[a,g(c)]$. Then $f$ is more repellent from
$a$ than $g$ on the set $B=\psi_g([g^2(c), g(c)])$ and so
$I_f\supset I_g$.

\endroster

\endproclaim

\demo{Proof} (1) Note that by the arguments we presented before
stating Lemma 3.1 $a$ is an $f$-fixed point.  Let $x$ be a periodic
point of $g$ of rotation pair $(p,q)$. Consider a $g$-loop
$[\psi_g(x),a], [\psi_g(g(x)),a],\dots, [\psi_g(g^{q-1}(x)),a]$ and
prove that it is an $f$-loop too. To this end it suffices to show
that $[f(\psi_g(g^i(x))),a]\supset [\psi_g(g^{i+1}(x)),a]$, and
indeed it follows immediately from the properties of $f$ and $g$ and
the construction of the function $\psi_g$.  Therefore by Lemma BM3 the
map $f$ also has a periodic point of rotation number $p/q$ which
proves that $I_f\supset I_g$ as desired. This easily implies the second part of
the statement (1) of the lemma.

(2) Clearly $\psi_g$ is the identity map on $[0,a')\cup [c,1]$. On
the other hand for any $x\in [a',c)$ we have $\psi_g(x)=x'_g$ where
$x'_g$ is a unique point such that $g(x'_g)=g(x), x'_g\neq x$. Hence
$B=\psi_g([g^2(c), g(c)])=[g^2(c), a']\cup [c, g(c)]$ and clearly
$f$ is more repellent from $a$ on $B$ than $g$. Due to the fact that
all periodic points of $g$ but the point $0$ are contained in a
$g$-invariant interval $[g^2(c), g(c)]$ it completes the proof.
{\unskip\nobreak\quad$\square$} \enddemo

We apply Lemma 3.1 to some pairs of maps and eventually to
one-parameter families of unimodal maps. Usually considered families
satisfy the conditions of pointwise growth or even are formed by
multiples of the same map; in other words the usual assumptions
differ from those of Lemma 3.1. However ``changing coordinates''
(i.e properly conjugating one of the maps involved) one can try to
see if two maps from a given family fit into the situation of Lemma
3.1.  We work with a few classes of maps having a single turning
point which for the sake of the definiteness is assumed to be
maximum from now on. Also without loss of generality we work from
now on with interval maps which map $0$ into itself and $1$ into
$0$.

Let $\Cal S=\{f: f$ is a convex (concave down) map of the interval $[0,1]$
into itself with a unique turning point $c_f$ which is maximum such
that $f|[0,c_f]$ and $f|[c_f,1]$ are $C^1$-maps, $f(c_f)>c_f$ and
$f(0)=f(1)=0 \}$ (for convenience we repeat here the definition given in
Introduction).  We can assume $f(c_f)>c_f$ without loss of
generality; under this assumption by the convexity of $f$ there is
no fixed point in $(0,c_f]$ and there is a unique fixed point
$a_f\in (c_f,1]$. For any $x\neq c_f$ there is a well defined point
$x'_f\neq x$ such that $f(x)=f(x'_f)$; also let $c'_f=c_f$.

\proclaim {Lemma 3.2} The following statements are true. \roster
\item Let $f,g\in \Cal S$ be such maps that $|f'|\ge |g'|$ and
$\dfrac {|c_f-a_f|}{|c_f-a'_f|}\ge \dfrac {|c_g-a_g|}{|c_g-a'_g|}$.
Then $I_f\supset I_g$.
\item Let $g\in \Cal S$ and $|g'(x)(x-c_g)|\le |g'(x'_g)(x'_g-c_g)|$
for any $x\ge c_g$. Then $I_{\nu g}\supset I_g$ for any $\nu>1$.
\endroster
\endproclaim

\demo{Proof} (1) We can assume that the set of periodic points
$P(g)$ of $g$ is not $\{1\}$. Observe that $c_g=c_f=c$ which follows
from $|f'|\ge |g'|$ since $0=|f'(c_f)|\ge |g'(c_f)|$. Notice also
that since $|f'|\ge |g'|$ then $a_g\le a_f$. Indeed, $|f'|\ge |g'|$
implies $\lambda (g[a_f,1])\le \lambda (f[a_f,1])$ (recall that by
$\lambda(A)$ we denote the Lebesgue measure of a set $A$); since
$f[a_f,1]=[0,a_f]$ and $g(1)=0$ we conclude that $g[a_f,1]\subset
[0,a_f]$ and so $g(a_f)\le a_f$ which implies that $a_g\le a_f$. Let
$\phi$ be a {\it linear} contraction of $[0,1]$ towards the point
$c$ onto its image $I'$ which maps $c$ into $c$ and $a_f$ into
$a_g=a$. The coefficient of contraction then is $q=\dfrac
{|c-a_g|}{|c-a_f|}\le 1$ so that $\phi'(u)=q$ for each $y$. The map
$\phi$ conjugates $f$ and $h:I'\to I'$; moreover, by construction
$g(a)=h(a)=a$. Let us show that $|h'(x)|\ge |g'(x)|$ for any $x\in
I'$. Indeed, points are attracted by $\phi$ closer to $c$ but remain
to the same side of $c$. Then by the properties of $f$ and $g$ and
by the choice of $\phi$ we have
$$|h'(x)|=q\cdot |f'(\phi^{-1}(x))|\cdot q^{-1}=|f'(\phi^{-1}(x))|\ge |g'(\phi^{-1}(x))| \ge |g'(x)|$$

Let us show that $h$ is more repellent from $a$ on $[c, h(c)]$ than
$g$.  Indeed, let $a\le y\le h(c)$. Then $\lambda (h[a,y])=a-h(y)\ge
\lambda (g[a,y])=a-g(y)$ which implies that $h(y)\le g(y)$.
Similarly for $c\le y\le a$ one can show that $h(y)\ge g(y)$.  In
particular it implies that $[g^2(c), g(c)]\subset [h^2(c), h(c)]$.
Let us check that $h$ is more repellent from $a$ than $g$ on
$[g^2(c), a'_g]$. First we show that $a'_g=a'\le a'_h$.  Indeed,
using the fact that $c=c_h=c_f=c_g$ and $a=a_g=a_h$ we can write
$\dfrac {|c-a|}{|c-a'_h|}=\dfrac {|c-a_f|}{|c-a'_f|}\ge \dfrac
{|c-a|}{|c-a'|}$ (the last inequality is actually given in the
conditions of the lemma as $a=a_g$ and $a'=a'_g$) which implies that
indeed $a'\le a'_h$. Hence $h(a')\le h(a'_h)=a=g(a')$. Since
$|h'|\ge |g'|$ we have $\lambda (h[z,a'])\ge \lambda (g[z,a'])$ for
any $z\in [g^2(c), a']$ which in turn implies that $h(z)\le g(z)$.
This allows to apply Lemma 3.1 and completes the proof of statement
(1). Figure 3.1 illustrates the proof. On this figure as well as on
some of the following figures the graphs of functions are given in
continuous or dashed lines while additional dotted lines are playing
explanatory role.

\define\figfour{\includegraphics{unim4.ps}}
\midinsert \moveright 0.375in \vbox{\hbox
{\hfil\lower3.7in\vbox{\figfour}\hfil}}
\botcaption {Figure 3.1} Functions $g$ (continuous line) and $f,h$
(dashed line)
\endcaption \endinsert

(2) This statement is in fact a corollary of statement (1). Indeed,
let $\nu g=f, c_g=c_f=c, a_g=a, a'_g=a'$; also for any $z$ let
$z'=z'_g=z'_f=\tau (z)$.  To apply statement (1) we need to check
that $\dfrac {|c-a_f|}{|c-a'_f|}\ge \dfrac {|c-a|}{|c-a'|}$. By the
proven in the beginning of the proof of (1), $a_f>a$. Thus it is
enough to show that for any pair of points $x,y$ such that $c<z\le
y$ we have $\dfrac {|c-y|}{|c-y'|}\ge \dfrac {|c-z|}{|c-z'|}$. Let
us show that the function $\psi(x)=\ln (x-c)/(c-x')=\ln (x-c) - \ln
(c-\tau(x))$ is increasing for $z\in (c,1]$.  Notice that
$\tau'(x)=g'(x)/g'(x')$; therefore $\psi'(x)=\dfrac 1{x-c} - \dfrac
{-\tau'(x)}{c-\tau(x)}= \dfrac 1{x-c} + \dfrac
{g'(x)}{g'(x')(c-x')}\ge 0$ (the last inequality follows from
$|g'(x)(x-c)|\le |g'(x')(x'-c)|$). {\unskip\nobreak\quad$\square$}
\enddemo

Call a map $g$ {\it even} if $g(x)=g(1-x)$ for any $x\in [0,1]$.
Clearly Lemma 3.2(1) holds if both $f$ and $g$ are even and $|f'|\ge
|g'|$ because then $\dfrac {|c_f-a_f|}{|c_f-a'_f|}=\dfrac
{|c_g-a_g|}{|c_g-a'_g|}=1$ automatically. E.g., this applies if $f$
and $g$ are from the quadratic family. Similarly, Lemma 3.2(2) holds
if $g$ is even.

For the sake of convenience let us call all polynomials of degree no
more than 3 ``cubic''. Another result close to that of Lemma 3.2
deals with maps from the class $\Cal G\subset \Cal S, \Cal G= \{g:
g\in \Cal S$ is a cubic polynomial map of $[0,1]$ into itself $ \}$.
Our purpose is to prove an analog of Lemma 3.2 for maps from $\Cal
G$. It is worth mentioning here that in fact the assumption of
convexity could be somewhat weakened; some statements are proven for
maps from the family $\Cal H=\{f: f$ is a cubic map of $[0,1]$ into
$\Bbb R_+ \cup 0$ with a unique critical point in $[0, 1]$ such that
$f(0)=f(1)=0 \}$. We combine geometrical trivia about cubic maps in
Lemma 3.3; also, some of the statements in what follows are proven
for the sake of completeness.

\proclaim { Lemma 3.3} Let $f\neq g$ be non-trivial cubic
polynomials. \roster
\item The graphs of $f$ and $g$
have no more than three common points and they cannot have two
common points at which their derivatives are equal.
\item If there are three different points $x<y<z$ at which $f=g$
then at none of them $f'=g'$ and $f-g$ has one sign on
$(-\infty,x)\cup (y,z)$ and the other on $(x,y)\cup (z, \infty)$.
\item If there are exactly two points $x<y$ between, say, $a$ and $b$
at which $f=g$ then either {\rm(a)} $f'(x)\neq g'(x), f'(y)\neq
g'(y)$ and $f-g$ has one  sign on $(a,x)\cup (y,b)$ and the other on
$(x,y)$, or {\rm(b)} $f'(x)=g'(x), f'(y)\neq g'(y)$ and $f-g$ has
one sign on $(a,x)\cup (x,y)$ and the other on $(y,b)$, or {\rm(c)}
$f'(x)\neq g'(x), f'(y)=g'(y)$ and $f-g$ has one sign on $(a,x)$ and
the other on $(x, y)\cup (y,b)$.
\item Let $a<b, f(a)<g(a), f(b)<g(b)$. Then either
{\rm (a)} $f(x)<g(x)$ for any $x\in (a,b)$, or {\rm (b)} there is a
single $x\in (a,b)$ such that $f(x)=g(x)$, $f(y)<g(y)$ for $y\in
(a,b)\setminus \{x\}$ and in fact $f'(x)=g'(x)$, or {\rm (c)} there
are two points $u<v$ in $(a,b)$ such that $f(y)<g(y)$ for $y\in
(a,u)$, $f(y)>g(y)$ for $y\in (u,v)$, $f(y)<g(y)$ for $y\in (v,b)$
and $f'(u)\neq g'(u), f'(v)\neq g'(v)$.
\item If $f\ge g$ on $[a,b]$ and $f(a)=g(a), f(b)=g(b)$ then in fact
$f>g$ on $(a,b)$.
\endroster
\endproclaim

\demo{Proof} Statements (1), (2), (3) follow immediately from the
fact that $f$ and $g$ are cubic. Statement (4) follows from (2) and
(3). Statement (5) follows from (2). {\unskip\nobreak\quad$\square$}
\enddemo

    The next lemma studies properties of conjugacies of cubic
maps.

\proclaim {Lemma 3.4} Let $f\ge g, f\neq g$ be two maps from $\Cal
H$. Let $\psi$ be a linear non-strict contraction with the fixed
point $0$ which maps $[0,1]$ onto $[0,v], v\le 1$ and conjugates $f$
to a map $h:[0,v]\to \Bbb R_+ \cup 0$. Then the following holds:
\roster
\item if\,\,\,$h\le g$ in a small right semi-neighborhood of
$0$ then $h'(0)=f'(0)=g'(0)$ and either $v=1, f=g=h$, or $h<g$ on
$(0, v]$;
\item if {\rm (1)} does not hold then
either {\rm(a)} $v=1$ and $h=f>g$ inside $(0,1)$, or
{\rm(b)} $v<1$ and there exists a point $u\in (0,v)$ such that $h(x)>g(x)$ if
$x\in (0,u)$, $h(x)<g(x)$ if $x\in (u,v]$ and $h'(u)\le g'(u)$.
\endroster

\endproclaim

\demo{Proof} Note that by the construction $h(v)=0<g(v)$ if $v<1$.

(1) By the assumption $h'(0)\le g'(0)$. On the other hand $f\ge g$
and so $f'(0)\ge g'(0)$. Since $f$ is conjugate to $h$ by a linear
map $\psi$ we see that $f'(0)=h'(0)$ and so $h'(0)\ge g'(0)$.  Thus
in fact $h'(0)=g'(0)=f'(0)$. If now $v=1$ we immediately get
$f=g=h$. Suppose that $v<1$. Then $h(v)=0<g(v)$. Let us show that
then $h<g$ on $(0, v]$. Indeed, otherwise either there are two
points between $0$ and $v$ at which $g=h$ (as $h'(0)=g'(0)$ and
$h(0)=g(0)$, this is impossible for cubic maps), or there is one
such point between $0$ and  $v$ and at that point both $g=h$ and
$g'=h'$ (which is impossible for similar reasons). Hence $h<g$ on
$(0, v]$.

(2) First assume that $v=1$ and so $h=f$; then by Lemma 3.3(5) we
have possibility (2)(a). So let $v<1$ and consider possibilities
concerning common points of the graphs of $h$ and $g$ on $(0,v]$. If
there are no such points then $h(v)<g(v)$ implies $h(x)<g(x)$ for
all $x\in (0,v]$ and by (1) we get possibility (2)(b) from the
lemma. Suppose there are such points. Then by Lemma 3.3(1) there are
no more than two of them. Consider a few cases.

(i) Let $0<s<t<v$ be such that $h(s)=g(s)$ and $h(t)=g(t)$.  By
Lemma 3.3.(2) then the fact that $h(v)<g(v)$ implies that
$h(x)<g(x)$ if $x\in (0,s)$. Hence by (1) $h'(0)=g'(0)$ which
contradicts Lemma 3.3(2).

(ii) Let $u\in (0,v]$ be the only point in this interval such that
$h(u)=g(u)$. Then $h<g$ on $(u, v]$. If $h<g$ on $(0, u)$ then
clearly $h'(u)=g'(u)$ and by (1) $h'(0)=g'(0)$. Since $g$ and $h$
are both cubic but not identical, this is impossible. Hence $h>g$ on
$(0, u)$. It follows then that in fact $h'(u)\le g'(u)$ and we get
possibility (b) from the lemma. This completes the proof.
{\unskip\nobreak\quad$\square$}
\enddemo

    In the next lemma we study some geometrical properties of maps
from $\Cal H$.

\proclaim {Lemma 3.5} Let $f\ge g, f\neq g$ be two maps from $\Cal
H$. Then the following holds. \roster
\item Let the line $y=\gamma x, \gamma>0$ intersect the graph of $g$ at a point
$(x, \gamma x), x>0$ and the graph of $f$ at a point $(z, \gamma z), z>x$.
Moreover, suppose that $g'(x)\le 0$. Then $f'(z)\le g'(x)$ and
moreover $x/x'_g<z/z'_f$.
\item $g(c_g)/c_g<f(c_f)/c_f$
\endroster

\endproclaim

\demo{Proof} (1) Consider the linear contraction $\psi$ with a fixed
point $0$ which maps $z$ into $x$ and conjugates $f$ and a map $h$.
Let $v=\psi(1)<1$. Then $h(x)=g(x)$ so neither Lemma 3.4(1) nor
Lemma 3.4(2)(a) holds. Thus, Lemma 3.4(2)(b) holds which implies
that $g>h$ on $(x, v)$ and $g<h$ on $(0, x)$. Moreover, by the same
lemma $h'(x)=f'(z)\le g'(x)$. The last inequality implies that,
since $g'(x)\le 0$ (i.e. $x\ge c_g$), then $f'(z)\le 0$ (i.e. $z\ge
c_f$).

Let $x'_g=x', x'_h=x'', z'_f=z'$. Then $x'\le x, x''<x$. By
construction $g(x')=g(x)=h(x)=h(x'')=\gamma x$. Since $g<h$ on $(0,
x)$ then $g(x'')<h(x'')=g(x)$ which implies $x''<x'$ since otherwise
$x'\le x''<x$ and so $g(x'')\ge g(x)$. Hence
$x/x'=x/x'_g<x/x''=x/x'_h=z/z'_f$.  Figure 3.2 illustrates the
arguments.

\define\figfive{\includegraphics{unim5.ps}}
\midinsert \moveright 0.375in \vbox{\hbox
{\hfil\lower3.9in\vbox{\figfive}\hfil}}
\botcaption {Figure 3.2} Functions $g$ (bold line), $f$ (regular
line) and $h$ (dashed line) from Lemma 3.5.
\endcaption \endinsert

(2) Consider the line $\alpha$ which connects the points $(0,0)$ and
$(c_g, g(c_g))$. Then the point $(c_g, f(c_g))$ is strictly above
this line. This implies that the graph of $f$ and the line $\alpha$
intersect strictly to the right of $c_g$. Suppose that the point of
intersection is $(z,y)$. Let us apply statement (1) to this
situation. Then $f'(z)\le g'(c_g)=0$ and so $z\ge c_f$. Since
$f(c_f)\ge f(z)=y$ we see that the point $(c_f, f(c_f))$ lies to the
left and above the point $(z,y)$.  Clearly it implies that
$g(c_g)/c_g\le f(c_f)/c_f$. {\unskip\nobreak\quad$\square$} \enddemo

The next lemma is important for the proof of the analog of Lemma 3.2
for maps from $\Cal G$.

\proclaim {Lemma 3.6} Let $f,g\in \Cal H$ and $g$ have a fixed point
$a_g\in [c_g,1]$. If $f\ge g, f\neq g$ then there is a unique
$f$-fixed point $a_f\ge c_f$ and $a_f>a_g$, $a_g/a'_g<a_f/a'_f$ and
$f'(a_f)\le g'(a_g)$.  \endproclaim

\demo {Proof} By Lemma 3.3(5) $f(x)>g(x)$ if $x\in (0,1)$. Consider
the line $y=x$. Clearly this line intersects the graph of $f$ at a
point $(z, w)$ where $z>a_g$ is an $f$-fixed point. Then by Lemma
3.5 $f'(z)\le g'(a_g)\le 0$ and so $z\ge c_f$. On the other hand
obviously if there exists an $f$-fixed point $a_f\in [c_g,1]$ then
it is unique. Thus, $z=a_f$, and by the above $a_f>a_g$. 
Observe that $f$ moves $a_g$ strictly to the right
($f(a_g)>g(a_g)=a_g$). It remains to apply Lemma 3.5(1) to the line
$y=x$ and points $(a_g,a_g)$ on the graph of $g$ and $(a_f,a_f)$ on
the graph of $f$. {\unskip\nobreak\quad$\square$}
\enddemo

Now we are ready to prove the analog of Lemma 3.2 for maps from
$\Cal G$.

\proclaim { Lemma 3.7} Let $f\ge g$ and $f,g\in \Cal G$. Then
$I_f\supset I_g$.
\endproclaim

\demo{Proof} The idea of the proof is the same as in Lemma 3.2. We
may assume that $f\neq g$ and that there is a fixed point of $g$ in
$[c_g, 1]$. Then by Lemma 3.6 $a_f>a_g$. Consider a contraction
$\phi$ which has a fixed point $0$ and maps $a_f$ to $a_g=a$ (thus,
$\phi(x)=(a_g/a_f)x$).  Let $\phi$ conjugate $f$ to a map
$h:[0,\phi(1)]\to [0,\phi(1)]$; by Lemma 3.6
$\phi(a'_f)=a'_f(a_g/a_f)<a'_g$. Consider now the linear contraction
$\psi$ which has a fixed point $a$ and maps $\phi(a'_f)$ to $a'_g$.
Then $\zeta=\psi\circ \phi$ is a linear map such that
$\zeta(a_f)=a_g=a, \zeta(a'_f)=a'_g=a', \zeta([0,1])=[u,v]\subset
[0,1]$. Let $w:[u,v]\to [u,v]$ be a map to which $f$ is conjugate by
$\zeta$. We prove that Lemma 3.1 is applicable to the maps $w$ and
$g$ which will then imply the required. Figure 3.3 illustrates the
proof.

\define\figsix{\includegraphics{unim6.ps}}
\midinsert \moveright 0.375in \vbox{\hbox
{\hfil\lower3.9in\vbox{\figsix}\hfil}}
\botcaption {Figure 3.3} Functions $g$ (bold line), $f$ (regular
line) and $w$ (dashed line) from Lemma 3.7
\endcaption \endinsert

Indeed, the convexity of $g$ implies that $w(u)=u<g(u)$.  Now, the
segment of straight line $\bold k$ connecting points $(a,a)$ and
$(\phi(1),0)$ is located below the graph of $g$ on the interval
$[a,\phi(1)]$ because  $g$ is convex. At the same time it is clear
that the linear map $\psi$ conjugates $h$ and $w$; since $\psi$ is
linear then a point $(\psi (\phi(1)), w(\psi(\phi(1))))=(v,w(v))$
lies on $\bold k$ and thus is still below the graph of $g$ which
means that $w(v)<g(v)$. These are the only two times we rely upon
the convexity in the proof.  Compare $w$ and $g$ on $[u,v]$. By
Lemma 3.3(4) we have $w(x)<g(x)$ if $x\in [u,a')$, $w(x)>g(x)$ if
$x\in (a',a)$ and $w(x)<g(x)$ if $x\in (a,v]$. This allows to apply
Lemma 3.2 which completes the proof. {\unskip\nobreak\quad$\square$}
\enddemo

    For the sake of convenience we sum up the results of Lemmas
3.2 and 3.7 in Theorem 3.8 dealing with one-parameter families of
interval maps.

\proclaim { Theorem 3.8} Let $f_\nu, \nu \in [b,d]$ be a
one-parameter family of interval maps such that one of the following
properties holds. \roster
\item $f_\nu\in \Cal S$ for any $\nu$; also,
if $\nu>\mu$ then $|f'_\nu|\ge |f'_\mu|$ and $\dfrac
{|c_{f_\nu}-a_{f_\nu}|}{|c_{f_\nu}-a'_{f_\nu}|}\ge \dfrac
{|c_{f_\mu}-a_{f_\mu}|}{|c_{f_\mu}-a'_{f_\mu}|}$.
\item $f=f_b\in \Cal S,\,|f'(x)(x-c_f)|\le |f'(x'_f)(x'_f-c_f)|$
for any $x\ge c_f$ and $f_\nu=\nu f$.
\item $f_\nu \in \Cal G$ for any $\nu$ and $f_\nu \ge f_\mu$ if $\nu>\mu$.
\endroster

Then $I_{f\nu}\supset I_{f_\mu}$ if $\nu>\mu$; in particular, if
$f_b$ has an odd periodic point then $P(f_\nu)\supset P(f_\mu)$ if
$\nu>\mu$.
\endproclaim

\enddocument